\documentclass[12pt]{article}
\usepackage{amssymb}
\usepackage{pstricks,pstricks-add,pst-math,pst-xkey}

\newcommand{\ba}{{\mbox{\bf a}}}
\newcommand{\bb}{{\mbox{\bf b}}}
\newcommand{\bab}{{\mbox{\bf ab}}}

\newcommand{\bc}{{\mbox{\bf c}}}
\newcommand{\lbc}{{\mbox{\rm\scriptsize\bf c}}}
\newcommand{\bd}{{\mbox{\bf d}}}
\newcommand{\lbd}{{\mbox{\rm\scriptsize\bf d}}}
\newcommand{\bcd}{{\mbox{\bf cd}}}

\newcommand{\BR}{{\mbox{\bf R}}}

\newcommand{\BZ}{{\mbox{\bf Z}}}
\newcommand{\CB}{{\mathcal B}}
\newcommand{\dd}{d}

\newcommand{\hh}{\hat h}

\newcommand{\lPhi}{\stackrel{\leftarrow}{\Phi}}

\newcommand{\lver}{{\mbox{\rm\scriptsize vert}}}
\newcommand{\of}{\overline f}
\newcommand{\oh}{\overline h}

\newcommand{\oR}{\overline R}

\newcommand{\rPhi}{\stackrel{\rightarrow}{\Phi}}

\newcommand{\Sum}{\displaystyle\sum}
\newcommand{\ver}{{\mbox{\rm vert}}}

\newtheorem{thm}{Theorem}
\newtheorem{lemma}{Lemma}

\newtheorem{cor}{Corollary}

\begin{document}

\begin{sloppypar}

\title{Sweeping the $\bcd$-Index and the Toric $h$-Vector}
\author{Carl W. Lee\\Department of Mathematics\\University of Kentucky}
\date{April 6, 2009, Revised November 1, 2010}
\maketitle

\noindent{\bf Abstract:}  
We derive 
formulas for the $\bcd$-index and the toric $h$-vector of a convex
polytope $P$ from a sweeping by a hyperplane.  
These arise from
interpreting the corresponding 
$S$-shelling~\cite{sta94:fla} of the dual of $P$.
We describe
a partition of the faces of the complete truncation of $P$ to
reflect explicitly the nonnegativity of its $\bcd$-index and what its
components are counting.
One corollary is a quick
way to compute the toric $h$-vector directly from the $\bcd$-index
that turns out to be an immediate consequence of formulas in
\cite{bayehr00:tor}.
We also propose an ``extended toric'' $h$-vector that fully captures
the information in the flag $h$-vector.

\section{Introduction}

By sweeping a hyperplane across a simple convex $\dd$-polytope $P$, the
$h$-vector, $h(P^*)=(h_0,\ldots,h_\dd)$, of its dual can be
computed---the edges in $P$ are
oriented in the direction of the sweep and $h_i$ equals the number of
vertices of outdegree $i$.  Moreover, the nonempty faces of $P$ can be
partitioned to explicitly reflect the formula for the $h$-vector.
For a general convex polytope,
in place of the $h$-vector, one often considers the flag $f$-vector
and flag $h$-vector as
well their encoding into the $\bcd$-index, and also the toric $h$-vector
(which does not contain the full information of the flag $h$-vector).
In this paper we derive 
formulas for the $\bcd$-index and for the toric $h$-vector
of a convex polytope $P$ from a sweeping of
$P$ (Theorems~\ref{cdsweep}, \ref{cdsymm},
\ref{cdtoh0} and \ref{cdtoh2}).
These arise from
interpreting the corresponding 
$S$-shelling~\cite{sta94:fla} of the dual of $P$.
We describe
a partition of the faces of the complete truncation of $P$ 
to
provide an interpretation of what 
the components of the $\bcd$-index are counting
(Theorem~\ref{partition} and Corollary~\ref{block}).
One corollary (Theorem~\ref{cdtoh1}) is a quick
way to compute the toric $h$-vector directly from the $\bcd$-index
that turns out to be an immediate consequence of formulas in
\cite{bayehr00:tor}.
We also propose an ``extended toric'' $h$-vector that fully captures
the information in the flag $h$-vector (Section~\ref{exttoric}).

Refer to
\cite{baylee93:com,bro83:int,gru03:con,hib92:alg,klekle95:con,mcmshe71:con,zie07:lec}, 
for example,
for background information on polytopes and their face numbers.

\section{The $h$-Vector}
\label{ordinaryh}

We begin by reviewing some well-known 
facts about $f$-vectors of polytopes.
For a convex $\dd$-dimensional polytope ($\dd$-polytope) $P$ let
$f_i=f_i(P)$ denote the number of $i$-faces ($i$-dimensional
faces) of $P$, $i=-1,\ldots,\dd$.  (Note that $f_{-1}=1$, counting the
empty set, and $f_\dd=1$, counting $P$ itself\@.)
The vector $f(P)=(f_0,\ldots,f_{\dd-1})$ is the {\em $f$-vector\/} of
$P$, and $f(P,x)$ is defined to be $\sum_{i=0}^\dd f_ix^i$.  
Faces of dimension $0$, $1$, and $\dd-1$ are called,
respectively, {\em
vertices}, {\em edges}, and {\em facets\/} of $P$.
The set of vertices of $P$ will be denoted $\ver(P)$.
A $\dd$-polytope is {\em simplicial\/} if every face
is a simplex.  
A $\dd$-polytope is {\em simple\/} if every vertex is contained in
exactly $\dd$ edges.  A dual to a simplicial polytope is simple, and
vice versa.

Let $P\subset\BR^\dd$ 
be a simple $\dd$-polytope.
The {\em $h$-vector\/} of the dual $P^*$ of $P$ is $(h_0,\ldots,h_d)$
where $h(P,x)=f(P,x-1)=\sum_{i=0}^dh_ix^i$.
Choose a direction $p\in\BR^\dd$ such that the inner product
$p\cdot x$ is different for each
vertex $v$ of $P$.  Sweep the hyperplane $H=\{x\in\BR^\dd:p\cdot x=q\}$
across $P$ by letting the parameter $q$ range from $-\infty$ to
$\infty$.  
(Recall that if $P$ contains the origin in its interior, then
ordering the vertices of 
$P$ using a sweeping hyperplane corresponds to
ordering the facets of the polar dual $P^*$
using a line shelling induced by a line
through the origin\@.)
Orient each edge of $P$ in the direction of increasing value of 
$p\cdot x$.

Each face of $P$ will have a unique minimal vertex with respect to
this orientation.  
To each vertex $v$ associate the set $B_v$ of
nonempty faces having $v$ as the minimal vertex, and (with a small
abuse of notation) associate
the monomial $h_v=x^i$,
where $i$ is the outdegree of $v$.
Then ${\CB}=\{B_v:v\in\ver(P)\}$ is a partition of the nonempty
faces of $P$.
The faces in 
$B_v$ contribute $(x+1)^i$ to $f(P,x)$ and so contribute $h_v$ to $h(P,x)$.
Therefore $h(P,x)=\sum_vh_v$ and each block $B_v$ contributes a
coefficient of $1$ to a single monomial.


%


\section{The $\bcd$-Index}

Two objects of study that each, in its own way, generalizes the simplicial
$h$-vector, are the $\bcd$-index and the toric $h$-vector.
Stanley~\cite{sta94:fla} introduced the notion of $S$-shellings to demonstrate
the nonnegativity of the $\bcd$-index.  

We will consider a sweeping of 
a polytope $P$ and, motivated by the
calculations associated with the $S$-shelling of its dual,
will construct a partition $\CB(P)$
of the nonempty faces of the complete
truncation of $P$, such that each block of $\CB(P)$ contributes a
coefficient of $1$ to
one word of the $\bcd$-index of $P$.  

\subsection{Definitions}

Let $P$ be a convex $\dd$-polytope.  Using the notation
$[\dd-1]=\{0,\ldots,\dd-1\}$, for every subset
$S=\{s_1,\ldots,s_k\}\subseteq[\dd-1]$ where 
$s_1<\cdots<s_k$,
define an {\em $S$-chain\/} to be a chain
of faces of $P$ of the form $F_1\subset\cdots\subset F_k$ where $F_i$
is face of $P$ of dimension $s_i$, $i=1,\ldots,k$.  Let $f_S(P)$ be the
number of $S$-chains.  The vector $\of(P)=(f_S(P))_{S\subseteq[\dd-1]}$ 
is the {\em flag $f$-vector\/} of $P$.

Now define 
\begin{equation}
\label{flagh}
h_S=h_S(P)=\sum_{T\subseteq S}(-1)^{|S|-|T|}f_T(P),\ S\subseteq[\dd-1].
\end{equation}
The vector $\oh(P)=(h_S(P))_{S\subseteq[\dd-1]}$ is the 
{\em flag $h$-vector\/} or {\em extended $h$-vector\/}
of $P$, introduced by Stanley~\cite{sta79:bal}.

Bayer and Billera showed that the affine span of the set
$\{\oh(P):h$ is a convex $\dd$-polytope$\}$ has dimension $F_\dd-1$,
where $F_\dd$ is 
the $\dd$th Fibonacci number.  Bayer and Klapper~\cite{baykla91:new} proved that
the flag $h$-vector can be encoded into the $\bcd$-index, which
precisely reflects this dimension.  Associate with each subset
$S\subseteq[\dd-1]$ the word $w_S=w_0\cdots w_{\dd-1}$ in the
noncommuting indeterminates $\ba$ and $\bb$, where $w_i=\ba$ if $i\not\in S$
and $w_i=\bb$ if $i\in S$.  The {\em $\bab$-index\/} of $P$ is then the
polynomial
\[
\Psi(P)=\Psi(P,\ba,\bb)=\sum_{S\subseteq[\dd-1]}h_S(P)w_S.
\]
The existence of the $\bcd$-index asserts that there is a polynomial in
the noncommuting indeterminates $\bc$ and $\bd$, 
$\Phi(P)=\Phi(P,\bc,\bd)$, such that
setting $\bc=\ba+\bb$ and $\bd=\bab+\bb\ba$ we have
$\Phi(P,\bc,\bd)=\Phi(P,\ba+\bb,\bab+\bb\ba)=\Psi(P,\ba,\bb)$.  
Note that $\bc$ has
degree one, $\bd$ has degree two,
and $\Phi(P)$ has degree $\dd$.  There are $F_\dd$ $\bcd$-words of degree
$\dd$ and one of them, $\bc^\dd$, will always have coefficient $1$.
Therefore the remaining $F_\dd-1$ terms of the $\bcd$-index capture the
dimension of the affine span of the flag $f$-vectors of $\dd$-polytopes.

%

\subsection{Partitioning the Complete Truncation}

Given a $d$-polytope, we will first construct its complete truncation $T(P)$, the faces of which are in bijection with the chains of $P$.  We will partition the faces of 
$T(P)$ into blocks, with a certain collection of blocks (and corresponding contribution toward $\Phi(P)$) associated with each vertex of $P$.  

Truncate all of the faces of $P$ by first
truncating the vertices of $P$, translating a supporting hyperplane
to each vertex a depth $\epsilon$ into $P$ and
giving each resulting $(\dd-1)$-face the label $0$.
Then continue by truncating
the original edges of $P$ at a depth of $\epsilon^2$ and
giving each resulting $(\dd-1)$-face the label $1$, truncating the
original $2$-faces of $P$ at a depth of $\epsilon^3$, etc., until 
finally truncating the
original $(d-1)$-faces of $P$ at a depth of $\epsilon^d$.  
Here, $\epsilon>0$ is taken to be sufficiently small for the sake of subsequent arguments.  
The resulting simple polytope, $T(P)$, called the
{\em complete truncation\/} of $P$,
is dual to the complete barycentric subdivision of the dual $P^*$ of $P$,
and its faces are in one-to-one correspondence with the chains of $P$. 
In fact, each nonempty face $G$ of $T(P)$ corresponds to an $S$-chain of $P$,
where $\sigma(G)=S$ is the set of labels of all of the
facets of $T(P)$ containing $G$.  The polytope $T(P)$ itself is 
labeled by
the empty set.  
For the sweeping hyperplane, choose a vector $p\in\BR^d$ such that the 
inner product 
$p\cdot x$ is different for all vertices occurring at all stages in the
truncation process.
See the first row of Figure~\ref{polygonpartition} for an example of a
pentagon and its truncation.

For each nonempty face $G$ of $T(P)$ of positive dimension $\dim(G)$ let
$j=\min\{i:i\not\in\sigma(G)\}$ and $w$ be the vertex
of $G$ with greatest value of $p\cdot x$.
Define the {\em top face\/} of $G$ to be  
the unique face $\tau(G)$ of $G$ of dimension $\dim(G)-1$
that contains $w$ and has 
label set $\sigma(G)\cup\{j\}$.
Similarly, let $w'$ be the vertex of $G$ with the smallest value of
$p\cdot x$, and define the {\em bottom face\/} of $G$ to be the unique face $\beta(G)$ of $G$ of dimension $\dim(G)-1$ that contains $w'$ and has 
label set $\sigma(G)\cup\{j\}$.
See the second row of Figure~\ref{polygonpartition}---each polygon
depicts a certain face of $T(P)$, together with its top and bottom
faces.

For vertex $v$ of $P$, let $Q_v$ be the $(d-1)$-face created when truncating $v$ in $P$, and $T(Q_v)$ be the complete truncation of $Q_v$ induced by $T(P)$.
Define $H_v=\{x\in\BR^d:p\cdot x=q_v\}$ 
to be the hyperplane in the sweeping family that contains
$v$, $H^+_v$ to be the open halfspace $\{x\in\BR^d:p\cdot x> q_v\}$, and 
$H^-_v$ to be the open halfspace $\{x\in\BR^d:p\cdot x<q_v\}$.
Faces of $T(Q_v)$ will be called {\em upper}, {\em middle}, or {\em
lower\/} faces according to whether they lie in $H^+_v$, intersect
$H_v$, or lie in $H^-_v$, respectively.
Note that if $v$ is the vertex of $P$ minimizing $p\cdot x$ then $T(Q_v)$ has no middle or lower faces, and if $v$ is the vertex of $P$ maximizing $p\cdot x$ then $T(Q_v)$ has no middle or upper faces.
Let $R_v$ be the polytope $Q_v\cap H_v$, which has dimension $d-2$ when it is nonempty.
($R_v$ will be empty if and only if $v$ minimizes or maximizes $p\cdot
x$ over $P$\@.)
Let $T(R_v)$ be the complete truncation of $R_v$ induced by $T(P)$; namely, $T(R_v)=T(Q_v)\cap H_v$.
Hence the faces of $T(R_v)$ are precisely the intersections of $H_v$ with the middle faces of $T(Q_v)$.
Observe that for a face $G$ of $T(P)$, $0\in\sigma(G)$ if and only if $G$ is a
face of some $T(Q_v)$.

\begin{lemma}
\label{nonzero}
For any face $G$ of $T(P)$ such that $0\not\in\sigma(G)$, the top face
$\tau(G)$ is a lower face of some $T(Q_v)$, and the bottom face
$\beta(G)$ is an upper face of some (other) $T(Q_v)$.  
Further, for every $v$, every lower and upper face of $T(Q_v)$ is uniquely obtainable in this way.
\end{lemma}

\noindent{\bf Proof.}
Suppose $0\not\in\sigma(G)$.  Then
$\sigma(\tau(G))=\sigma(G)\cup\{0\}$.
Let $v$ be the vertex of $P$ for which $T(Q_v)$ contributes the label
$\{0\}$ to $\tau(G)$, and
let $w$ be the vertex of $G$ that maximizes $p\cdot x$ over $G$.
Then $p\cdot w<p\cdot v$, and so $\tau(G)$, which is a face of
$T(Q_v)$, lies in $H^-_v$.  The analogous argument shows that
$\beta(G)$ is an upper face of some $T(Q_v)$.
Now let $G'$ be a lower face of some $T(Q_v)$.
$G'$ corresponds to an $S$-chain $F_1\subset\cdots\subset F_k$ in $P$, 
$S=\{s_1,\ldots,s_k\}$, where 
$0=s_1<s_2<\cdots<s_k$ and 
$F_1=\{v\}$.
Each $F_i$ contributes a facet $F'_i$ to $T(P)$ and $G'$ is the intersection
of these facets.
Because $G'$ lies in $H^-_v$, by convexity we conclude that there is
some $F'_i\not=F'_1$ that also lies in $H^-_v$.
Define $G$ to be the unique face of $T(P)$ with label set
$\sigma(G)=\sigma(G')\setminus\{0\}$ that contains $G'$.  Then
$G=F'_2\cap\cdots\cap F'_k$ lies in $H^-_v$.
Hence the top vertex of $G$ cannot lie above $H_v$ or be associated
with any $T_{v'}$ for any higher vertex $v'$ of $P$, and so must be in
$G'$, confirming that $G'=\tau(G)$.~$\Box$

Given the partitions for complete truncations of polytopes of
dimension less than $d$, 
we will recursively define the partition $\CB(P)$ of the faces of $T(P)$.
Three properties to be maintained are:

\begin{description}
\item[P1.]
Every vertex $v$ of $P$ will contribute an associated (though possibly
empty) collection $\CB_v(P)$ of blocks to the partition.  
\item[P2.]
If $d>0$ then every face $G$ for which $0\not\in\sigma(G)$ will be
placed in the same block as its top face $\tau(G)$.
\item[P3.]
Suppose $d>0$ and $H$ is any hyperplane in the sweeping family not
meeting any $T(Q_v)$.
Then for any vertex $v$ of $P$ in
$H^+$, the faces in the blocks $\CB_v(P)$ all lie in $H^+$.
\end{description}

\noindent{\bf Construction of $\CB(P)$:}
\begin{description}
\item[Step 0:]  If $P$ is a $0$-polytope, $T(P)$ is a single vertex $v$ 
and $\CB_v(P)$ contains the single block $\{v\}$.
So assume that $P$ has positive dimension.  
\item[Step 1:]  For every face $G$ of $T(P)$ such that
$0\not\in\sigma(G)$ create the ``pre-block'' $\{G,\tau(G),\beta(G)\}$
consisting of $G$, its top face and its bottom face.  At this point,
by Lemma~\ref{nonzero}, every face of $T(P)$ except the middle faces of
the various $T(Q_v)$ have been assigned to pre-blocks.  
\item[Step 2:]  For each vertex $v$ and each middle face $G$ of
$T(Q_v)$, insert $G$ in the pre-block containing its top face $\tau(G)$, which will be an upper face of $T(Q_v)$.  
At this point every face of $T(P)$ has been assigned to a pre-block,
there is a one-to-one correspondence between upper faces and
pre-blocks, and middle faces are in separate pre-blocks.
\item[Step 3:]  For each vertex $v$, consider the recursively defined
partition $\CB(R_v)$
of the faces of $T(R_v)$ (empty if $R_v$ is empty).  Let $B$
be a block in this partition.  Each face in $B$ corresponds to a
certain middle face in $T(Q_v)$.  
Merge the pre-blocks containing these corresponding middle
faces into a block $B'$.  Place $B'$ into $\CB_v(P)$.
\item[Step 4:]  For each vertex $v$,
consider the recursively defined partition $\CB(Q_v)$ of the faces of
$T(Q_v)$.
For each vertex $w$ of $Q_v$ in $H^+_v$, let $\CB_w(Q_v)$ be the
blocks of $\CB(Q_v)$ associated with $w$.
Let $B$ be a block in $\CB_w(Q_v)$ (if any).  By property (P3)
the faces in $B$ are certain
upper faces of $T(Q_v)$.  Merge the pre-blocks containing these upper
faces into a block $B'$, and place $B'$ into $\CB_v(P)$.
Once this is carried out for every vertex $v$ of $P$,
all of the pre-blocks have been merged as necessary and
$\CB(P)=\bigcup_v\CB_v(P)$.
\end{description}

To verify that there are no conflicts between the mergings in Step~3
and the mergings in Step~4, we need to make some observations.  
Let $G$ be a middle face of $T(Q_v)$.  Note that $0\in\sigma(G)$ but
$1\not\in\sigma(G)$, because $H_v$ does not contain any vertices of
$Q_v$ and the truncations of the edges and other faces of $P$ are made at 
sufficiently small depths.  Now regard $Q_v$ as a polytope in its own
right.  The label set $\sigma'(G)$ of $G$ with respect to $T(Q_v)$ is
obtained from that of $\sigma(G)$ by deleting $0$ and reducing the
remaining elements of $\sigma(G)$ by one.  Thus $0\not\in\sigma'(G)$.
By property (P2), within $\CB(Q_v)$, $G$ will be placed in the same
block as $\tau(G)$.  Thus the blocks in $\CB(P)$, restricted to the faces
in $Q_v$, will be blocks or subsets of blocks in the partition of
the faces of $T(Q_v)$.

It is straightforward from the construction 
to verify that $\CB(P)$ satisfies properties (P1)--(P3).

\begin{thm}
\label{partition}
$\CB(P)$ is a partition of $T(P)$.
\end{thm}

\noindent{\bf Examples}
\begin{enumerate}
\item The line segment ($d=1$).  See Figure~\ref{segment}.
\begin{figure}
\begin{center}
\psset{xunit=1.0cm,yunit=1.0cm,algebraic=true,dotstyle=*,dotsize=3pt
0,linewidth=0.8pt,arrowsize=3pt 2,arrowinset=0.25}
\begin{pspicture*}(-1,1.04)(3,4.64)
\psline(0,4)(0,2)
\psline[linewidth=5.2pt](2,3.68)(2,2.3)
\rput[tl](-0.74,3.9){$Q_{v_2}$}
\rput[tl](-0.76,2.56){$Q_{v_1}$}
\rput[tl](2.34,2.46){$\bc$}
\rput[tl](-0.1,4.58){$v_2$}
\rput[tl](-0.08,1.78){$v_1$}
\psdots(0,4)
\psdots(0,2)
\psdots(0,3.68)
\psdots(0,2.3)
\psdots[dotsize=9pt 0,dotstyle=*](2,3.68)
\psdots[dotsize=9pt 0,dotstyle=*](2,2.3)
\end{pspicture*}
\end{center}
\caption{Partitioning the Truncation of a Line Segment}
\label{segment}
\end{figure}
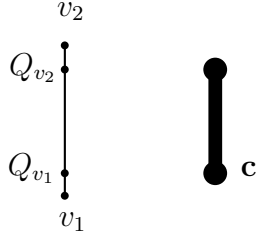
If $P$ is a line segment with two
vertices swept in the order $v_1,v_2$, then $Q_{v_i}$ is a point and
$R_{v_i}$ is empty, $i=1,2$.  There is only one pre-block, and this
becomes the only block in the partition of $T(P)$.
\item The $n$-gon ($\dd=2$).  See Figures~\ref{polygonpartition} and
\ref{polygon}.
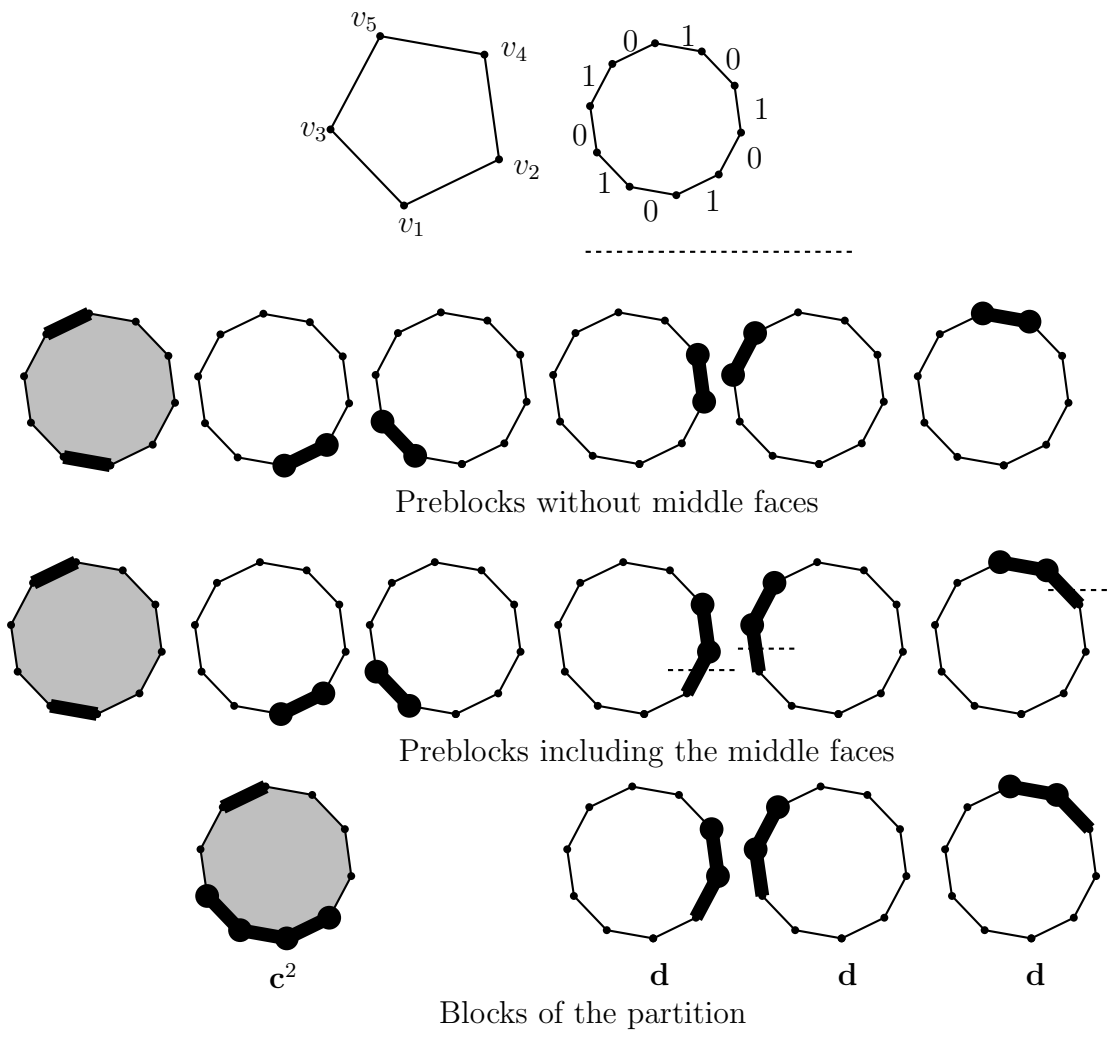
\begin{figure}
\begin{center}
\hspace*{-3cm}
\psset{xunit=.5cm,yunit=.5cm,algebraic=true,dotstyle=o,dotsize=3pt
0,linewidth=0.8pt,arrowsize=3pt 2,arrowinset=0.25}
\begin{pspicture*}(-15.95,-24.37)(34.33,5.2)
\pspolygon[linestyle=none,fillstyle=solid,fillcolor=lightgray](-6.85,-7.25)(-5.61,-7.47)(-4.48,-6.92)(-3.89,-5.81)(-4.06,-4.56)(-4.93,-3.65)(-6.17,-3.43)(-7.31,-3.98)(-7.9,-5.1)(-7.72,-6.34)
\pspolygon[linestyle=none,fillstyle=solid,fillcolor=lightgray](-2.16,-19.84)(-0.92,-20.06)(0.21,-19.51)(0.8,-18.4)(0.63,-17.15)(-0.24,-16.24)(-1.48,-16.02)(-2.62,-16.57)(-3.21,-17.69)(-3.03,-18.93)
\pspolygon[linestyle=none,fillstyle=solid,fillcolor=lightgray](-7.2,-13.87)(-5.96,-14.09)(-4.83,-13.54)(-4.24,-12.43)(-4.41,-11.18)(-5.28,-10.27)(-6.52,-10.05)(-7.66,-10.6)(-8.25,-11.72)(-8.07,-12.96)
\psline(1.57,3.95)(4.34,3.46)
\psline(4.34,3.46)(4.73,0.67)
\psline(4.73,0.67)(2.2,-0.56)
\psline(2.2,-0.56)(0.25,1.47)
\psline(1.57,3.95)(0.25,1.47)
\psline(8.2,-0.06)(9.44,-0.28)
\psline(9.44,-0.28)(10.57,0.27)
\psline(10.57,0.27)(11.16,1.38)
\psline(11.16,1.38)(10.99,2.63)
\psline(10.99,2.63)(10.12,3.54)
\psline(10.12,3.54)(8.88,3.76)
\psline(8.88,3.76)(7.74,3.21)
\psline(7.74,3.21)(7.15,2.09)
\psline(7.15,2.09)(7.33,0.85)
\psline(7.33,0.85)(8.2,-0.06)
\psline[linewidth=5.2pt](-6.85,-7.25)(-5.61,-7.47)
\psline(-5.61,-7.47)(-4.48,-6.92)
\psline(-4.48,-6.92)(-3.89,-5.81)
\psline(-3.89,-5.81)(-4.06,-4.56)
\psline(-4.06,-4.56)(-4.93,-3.65)
\psline(-4.93,-3.65)(-6.17,-3.43)
\psline[linewidth=5.2pt](-6.17,-3.43)(-7.31,-3.98)
\psline(-7.31,-3.98)(-7.9,-5.1)
\psline(-7.9,-5.1)(-7.72,-6.34)
\psline(-7.72,-6.34)(-6.85,-7.25)
\psline(-2.22,-7.26)(-0.98,-7.48)
\psline[linewidth=5.2pt](-0.98,-7.48)(0.15,-6.93)
\psline(0.15,-6.93)(0.74,-5.82)
\psline(0.74,-5.82)(0.57,-4.57)
\psline(0.57,-4.57)(-0.3,-3.66)
\psline(-0.3,-3.66)(-1.54,-3.44)
\psline(-1.54,-3.44)(-2.68,-3.99)
\psline(-2.68,-3.99)(-3.27,-5.11)
\psline(-3.27,-5.11)(-3.09,-6.35)
\psline(-3.09,-6.35)(-2.22,-7.26)
\psline(2.5,-7.22)(3.74,-7.44)
\psline(3.74,-7.44)(4.87,-6.89)
\psline(4.87,-6.89)(5.46,-5.78)
\psline(5.46,-5.78)(5.29,-4.53)
\psline(5.29,-4.53)(4.42,-3.62)
\psline(4.42,-3.62)(3.18,-3.4)
\psline(3.18,-3.4)(2.04,-3.95)
\psline(2.04,-3.95)(1.45,-5.07)
\psline(1.45,-5.07)(1.63,-6.31)
\psline[linewidth=5.2pt](1.63,-6.31)(2.5,-7.22)
\psline(7.22,-7.22)(8.46,-7.44)
\psline(8.46,-7.44)(9.59,-6.89)
\psline(9.59,-6.89)(10.18,-5.78)
\psline[linewidth=5.2pt](10.18,-5.78)(10.01,-4.53)
\psline(10.01,-4.53)(9.14,-3.62)
\psline(9.14,-3.62)(7.9,-3.4)
\psline(7.9,-3.4)(6.76,-3.95)
\psline(6.76,-3.95)(6.17,-5.07)
\psline(6.17,-5.07)(6.35,-6.31)
\psline(6.35,-6.31)(7.22,-7.22)
\psline(12,-7.22)(13.24,-7.44)
\psline(13.24,-7.44)(14.37,-6.89)
\psline(14.37,-6.89)(14.96,-5.78)
\psline(14.96,-5.78)(14.79,-4.53)
\psline(14.79,-4.53)(13.92,-3.62)
\psline(13.92,-3.62)(12.68,-3.4)
\psline(12.68,-3.4)(11.54,-3.95)
\psline[linewidth=5.2pt](11.54,-3.95)(10.95,-5.07)
\psline(10.95,-5.07)(11.13,-6.31)
\psline(11.13,-6.31)(12,-7.22)
\psline(16.91,-7.25)(18.15,-7.47)
\psline(18.15,-7.47)(19.28,-6.92)
\psline(19.28,-6.92)(19.87,-5.81)
\psline(19.87,-5.81)(19.7,-4.56)
\psline(19.7,-4.56)(18.83,-3.65)
\psline[linewidth=5.2pt](18.83,-3.65)(17.59,-3.43)
\psline(17.59,-3.43)(16.45,-3.98)
\psline(16.45,-3.98)(15.86,-5.1)
\psline(15.86,-5.1)(16.04,-6.34)
\psline(16.04,-6.34)(16.91,-7.25)
\psline[linewidth=5.2pt](-2.16,-19.84)(-0.92,-20.06)
\psline[linewidth=5.2pt](-0.92,-20.06)(0.21,-19.51)
\psline(0.21,-19.51)(0.8,-18.4)
\psline(0.8,-18.4)(0.63,-17.15)
\psline(0.63,-17.15)(-0.24,-16.24)
\psline(-0.24,-16.24)(-1.48,-16.02)
\psline[linewidth=5.2pt](-1.48,-16.02)(-2.62,-16.57)
\psline(-2.62,-16.57)(-3.21,-17.69)
\psline(-3.21,-17.69)(-3.03,-18.93)
\psline[linewidth=5.2pt](-3.03,-18.93)(-2.16,-19.84)
\psline(7.59,-19.84)(8.83,-20.06)
\psline(8.83,-20.06)(9.96,-19.51)
\psline[linewidth=5.2pt](9.96,-19.51)(10.55,-18.4)
\psline[linewidth=5.2pt](10.55,-18.4)(10.38,-17.15)
\psline(10.38,-17.15)(9.51,-16.24)
\psline(9.51,-16.24)(8.27,-16.02)
\psline(8.27,-16.02)(7.13,-16.57)
\psline(7.13,-16.57)(6.54,-17.69)
\psline(6.54,-17.69)(6.72,-18.93)
\psline(6.72,-18.93)(7.59,-19.84)
\psline(12.6,-19.84)(13.84,-20.06)
\psline(13.84,-20.06)(14.97,-19.51)
\psline(14.97,-19.51)(15.56,-18.4)
\psline(15.56,-18.4)(15.39,-17.15)
\psline(15.39,-17.15)(14.52,-16.24)
\psline(14.52,-16.24)(13.28,-16.02)
\psline(13.28,-16.02)(12.14,-16.57)
\psline[linewidth=5.2pt](12.14,-16.57)(11.55,-17.69)
\psline[linewidth=5.2pt](11.55,-17.69)(11.73,-18.93)
\psline(11.73,-18.93)(12.6,-19.84)
\psline(17.64,-19.84)(18.88,-20.06)
\psline(18.88,-20.06)(20.01,-19.51)
\psline(20.01,-19.51)(20.6,-18.4)
\psline(20.6,-18.4)(20.43,-17.15)
\psline[linewidth=5.2pt](20.43,-17.15)(19.56,-16.24)
\psline[linewidth=5.2pt](19.56,-16.24)(18.32,-16.02)
\psline(18.32,-16.02)(17.18,-16.57)
\psline(17.18,-16.57)(16.59,-17.69)
\psline(16.59,-17.69)(16.77,-18.93)
\psline(16.77,-18.93)(17.64,-19.84)
\psline[linewidth=5.2pt](-7.2,-13.87)(-5.96,-14.09)
\psline(-5.96,-14.09)(-4.83,-13.54)
\psline(-4.83,-13.54)(-4.24,-12.43)
\psline(-4.24,-12.43)(-4.41,-11.18)
\psline(-4.41,-11.18)(-5.28,-10.27)
\psline(-5.28,-10.27)(-6.52,-10.05)
\psline[linewidth=5.2pt](-6.52,-10.05)(-7.66,-10.6)
\psline(-7.66,-10.6)(-8.25,-11.72)
\psline(-8.25,-11.72)(-8.07,-12.96)
\psline(-8.07,-12.96)(-7.2,-13.87)
\psline(-2.31,-13.87)(-1.07,-14.09)
\psline[linewidth=5.2pt](-1.07,-14.09)(0.06,-13.54)
\psline(0.06,-13.54)(0.65,-12.43)
\psline(0.65,-12.43)(0.48,-11.18)
\psline(0.48,-11.18)(-0.39,-10.27)
\psline(-0.39,-10.27)(-1.63,-10.05)
\psline(-1.63,-10.05)(-2.77,-10.6)
\psline(-2.77,-10.6)(-3.36,-11.72)
\psline(-3.36,-11.72)(-3.18,-12.96)
\psline(-3.18,-12.96)(-2.31,-13.87)
\psline(2.34,-13.87)(3.58,-14.09)
\psline(3.58,-14.09)(4.71,-13.54)
\psline(4.71,-13.54)(5.3,-12.43)
\psline(5.3,-12.43)(5.13,-11.18)
\psline(5.13,-11.18)(4.26,-10.27)
\psline(4.26,-10.27)(3.02,-10.05)
\psline(3.02,-10.05)(1.88,-10.6)
\psline(1.88,-10.6)(1.29,-11.72)
\psline(1.29,-11.72)(1.47,-12.96)
\psline[linewidth=5.2pt](1.47,-12.96)(2.34,-13.87)
\psline(7.35,-13.87)(8.59,-14.09)
\psline(8.59,-14.09)(9.72,-13.54)
\psline[linewidth=5.2pt](9.72,-13.54)(10.31,-12.43)
\psline[linewidth=5.2pt](10.31,-12.43)(10.14,-11.18)
\psline(10.14,-11.18)(9.27,-10.27)
\psline(9.27,-10.27)(8.03,-10.05)
\psline(8.03,-10.05)(6.89,-10.6)
\psline(6.89,-10.6)(6.3,-11.72)
\psline(6.3,-11.72)(6.48,-12.96)
\psline(6.48,-12.96)(7.35,-13.87)
\psline(12.51,-13.87)(13.75,-14.09)
\psline(13.75,-14.09)(14.88,-13.54)
\psline(14.88,-13.54)(15.47,-12.43)
\psline(15.47,-12.43)(15.3,-11.18)
\psline(15.3,-11.18)(14.43,-10.27)
\psline(14.43,-10.27)(13.19,-10.05)
\psline(13.19,-10.05)(12.05,-10.6)
\psline[linewidth=5.2pt](12.05,-10.6)(11.46,-11.72)
\psline[linewidth=5.2pt](11.46,-11.72)(11.64,-12.96)
\psline(11.64,-12.96)(12.51,-13.87)
\psline(17.37,-13.87)(18.61,-14.09)
\psline(18.61,-14.09)(19.74,-13.54)
\psline(19.74,-13.54)(20.33,-12.43)
\psline(20.33,-12.43)(20.16,-11.18)
\psline[linewidth=5.2pt](20.16,-11.18)(19.29,-10.27)
\psline[linewidth=5.2pt](19.29,-10.27)(18.05,-10.05)
\psline(18.05,-10.05)(16.91,-10.6)
\psline(16.91,-10.6)(16.32,-11.72)
\psline(16.32,-11.72)(16.5,-12.96)
\psline(16.5,-12.96)(17.37,-13.87)
\rput[tl](5.11,0.64){$v_2$}
\rput[tl](-0.56,1.66){$v_3$}
\rput[tl](4.78,3.79){$v_4$}
\rput[tl](0.79,4.57){$v_5$}
\rput[tl](8.56,-0.44){0}
\rput[tl](10.21,-0.17){1}
\rput[tl](11.32,0.97){0}
\rput[tl](11.5,2.29){1}
\rput[tl](10.75,3.61){0}
\rput[tl](9.55,4.24){1}
\rput[tl](8.02,4.09){0}
\rput[tl](6.91,3.16){1}
\rput[tl](6.67,1.6){0}
\rput[tl](7.33,0.31){1}
\rput[tl](1.97,-8.12){Preblocks without middle faces}
\rput[tl](2.06,-14.75){Preblocks including the middle faces}
\rput[t1](-1,-20.72){$\bc^2$}
\rput[t1](9,-20.72){$\bd$}
\rput[t1](14,-20.72){$\bd$}
\rput[t1](19,-20.72){$\bd$}
\rput[tl](3.14,-21.72){Blocks of the partition}
\psline[linestyle=dashed,dash=2pt 2pt](7.03,-1.79)(14.11,-1.79)
\psline[linestyle=dashed,dash=2pt 2pt](9.22,-12.92)(10.99,-12.92)
\psline[linestyle=dashed,dash=2pt 2pt](11.08,-12.35)(12.58,-12.35)
\psline[linestyle=dashed,dash=2pt 2pt](19.33,-10.79)(20.89,-10.79)
\rput[tl](2.05,-0.89){$v_1$}
\psdots[dotstyle=*](1.57,3.95)
\psdots[dotstyle=*](4.34,3.46)
\psdots[dotstyle=*](4.73,0.67)
\psdots[dotstyle=*](2.2,-0.56)
\psdots[dotstyle=*](0.25,1.47)
\psdots[dotstyle=*](8.2,-0.06)
\psdots[dotstyle=*](9.44,-0.28)
\psdots[dotstyle=*](10.57,0.27)
\psdots[dotstyle=*](11.16,1.38)
\psdots[dotstyle=*](10.99,2.63)
\psdots[dotstyle=*](10.12,3.54)
\psdots[dotstyle=*](8.88,3.76)
\psdots[dotstyle=*](7.74,3.21)
\psdots[dotstyle=*](7.15,2.09)
\psdots[dotstyle=*](7.33,0.85)
\psdots[dotstyle=*](-5.61,-7.47)
\psdots[dotsize=9pt 0,dotstyle=*](-0.98,-7.48)
\psdots[dotstyle=*](3.74,-7.44)
\psdots[dotstyle=*](8.46,-7.44)
\psdots[dotstyle=*](13.24,-7.44)
\psdots[dotstyle=*](18.15,-7.47)
\psdots[dotstyle=*](-6.85,-7.25)
\psdots[dotstyle=*](-5.61,-7.47)
\psdots[dotstyle=*](-4.48,-6.92)
\psdots[dotstyle=*](-3.89,-5.81)
\psdots[dotstyle=*](-4.06,-4.56)
\psdots[dotstyle=*](-4.93,-3.65)
\psdots[dotstyle=*](-6.17,-3.43)
\psdots[dotstyle=*](-7.31,-3.98)
\psdots[dotstyle=*](-7.9,-5.1)
\psdots[dotstyle=*](-7.72,-6.34)
\psdots[dotstyle=*](-2.22,-7.26)
\psdots[dotstyle=*](-0.98,-7.48)
\psdots[dotsize=9pt 0,dotstyle=*](0.15,-6.93)
\psdots[dotstyle=*](0.74,-5.82)
\psdots[dotstyle=*](0.57,-4.57)
\psdots[dotstyle=*](-0.3,-3.66)
\psdots[dotstyle=*](-1.54,-3.44)
\psdots[dotstyle=*](-2.68,-3.99)
\psdots[dotstyle=*](-3.27,-5.11)
\psdots[dotstyle=*](-3.09,-6.35)
\psdots[dotsize=9pt 0,dotstyle=*](2.5,-7.22)
\psdots[dotstyle=*](3.74,-7.44)
\psdots[dotstyle=*](4.87,-6.89)
\psdots[dotstyle=*](5.46,-5.78)
\psdots[dotstyle=*](5.29,-4.53)
\psdots[dotstyle=*](4.42,-3.62)
\psdots[dotstyle=*](3.18,-3.4)
\psdots[dotstyle=*](2.04,-3.95)
\psdots[dotstyle=*](1.45,-5.07)
\psdots[dotsize=9pt 0,dotstyle=*](1.63,-6.31)
\psdots[dotstyle=*](7.22,-7.22)
\psdots[dotstyle=*](8.46,-7.44)
\psdots[dotstyle=*](9.59,-6.89)
\psdots[dotsize=9pt 0,dotstyle=*](10.18,-5.78)
\psdots[dotsize=9pt 0,dotstyle=*](10.01,-4.53)
\psdots[dotstyle=*](9.14,-3.62)
\psdots[dotstyle=*](7.9,-3.4)
\psdots[dotstyle=*](6.76,-3.95)
\psdots[dotstyle=*](6.17,-5.07)
\psdots[dotstyle=*](6.35,-6.31)
\psdots[dotstyle=*](12,-7.22)
\psdots[dotstyle=*](13.24,-7.44)
\psdots[dotstyle=*](14.37,-6.89)
\psdots[dotstyle=*](14.96,-5.78)
\psdots[dotstyle=*](14.79,-4.53)
\psdots[dotstyle=*](13.92,-3.62)
\psdots[dotstyle=*](12.68,-3.4)
\psdots[dotsize=9pt 0,dotstyle=*](11.54,-3.95)
\psdots[dotsize=9pt 0,dotstyle=*](10.95,-5.07)
\psdots[dotstyle=*](11.13,-6.31)
\psdots[dotstyle=*](16.91,-7.25)
\psdots[dotstyle=*](18.15,-7.47)
\psdots[dotstyle=*](19.28,-6.92)
\psdots[dotstyle=*](19.87,-5.81)
\psdots[dotstyle=*](19.7,-4.56)
\psdots[dotsize=9pt 0,dotstyle=*](18.83,-3.65)
\psdots[dotsize=9pt 0,dotstyle=*](17.59,-3.43)
\psdots[dotstyle=*](16.45,-3.98)
\psdots[dotstyle=*](15.86,-5.1)
\psdots[dotstyle=*](16.04,-6.34)
\psdots[dotsize=9pt 0,dotstyle=*](-0.92,-20.06)
\psdots[dotstyle=*](8.83,-20.06)
\psdots[dotstyle=*](13.84,-20.06)
\psdots[dotstyle=*](18.88,-20.06)
\psdots[dotsize=9pt 0,dotstyle=*](-2.16,-19.84)
\psdots[dotstyle=*](-0.92,-20.06)
\psdots[dotsize=9pt 0,dotstyle=*](0.21,-19.51)
\psdots[dotstyle=*](0.8,-18.4)
\psdots[dotstyle=*](0.63,-17.15)
\psdots[dotstyle=*](-0.24,-16.24)
\psdots[dotstyle=*](-1.48,-16.02)
\psdots[dotstyle=*](-2.62,-16.57)
\psdots[dotstyle=*](-3.21,-17.69)
\psdots[dotsize=9pt 0,dotstyle=*](-3.03,-18.93)
\psdots[dotstyle=*](7.59,-19.84)
\psdots[dotstyle=*](8.83,-20.06)
\psdots[dotstyle=*](9.96,-19.51)
\psdots[dotsize=9pt 0,dotstyle=*](10.55,-18.4)
\psdots[dotsize=9pt 0,dotstyle=*](10.38,-17.15)
\psdots[dotstyle=*](9.51,-16.24)
\psdots[dotstyle=*](8.27,-16.02)
\psdots[dotstyle=*](7.13,-16.57)
\psdots[dotstyle=*](6.54,-17.69)
\psdots[dotstyle=*](6.72,-18.93)
\psdots[dotstyle=*](12.6,-19.84)
\psdots[dotstyle=*](13.84,-20.06)
\psdots[dotstyle=*](14.97,-19.51)
\psdots[dotstyle=*](15.56,-18.4)
\psdots[dotstyle=*](15.39,-17.15)
\psdots[dotstyle=*](14.52,-16.24)
\psdots[dotstyle=*](13.28,-16.02)
\psdots[dotsize=9pt 0,dotstyle=*](12.14,-16.57)
\psdots[dotsize=9pt 0,dotstyle=*](11.55,-17.69)
\psdots[dotstyle=*](11.73,-18.93)
\psdots[dotstyle=*](17.64,-19.84)
\psdots[dotstyle=*](18.88,-20.06)
\psdots[dotstyle=*](20.01,-19.51)
\psdots[dotstyle=*](20.6,-18.4)
\psdots[dotstyle=*](20.43,-17.15)
\psdots[dotsize=9pt 0,dotstyle=*](19.56,-16.24)
\psdots[dotsize=9pt 0,dotstyle=*](18.32,-16.02)
\psdots[dotstyle=*](17.18,-16.57)
\psdots[dotstyle=*](16.59,-17.69)
\psdots[dotstyle=*](16.77,-18.93)
\psdots[dotstyle=*](-5.96,-14.09)
\psdots[dotsize=9pt 0,dotstyle=*](-1.07,-14.09)
\psdots[dotstyle=*](3.58,-14.09)
\psdots[dotstyle=*](8.59,-14.09)
\psdots[dotstyle=*](13.75,-14.09)
\psdots[dotstyle=*](18.61,-14.09)
\psdots[dotstyle=*](-7.2,-13.87)
\psdots[dotstyle=*](-5.96,-14.09)
\psdots[dotstyle=*](-4.83,-13.54)
\psdots[dotstyle=*](-4.24,-12.43)
\psdots[dotstyle=*](-4.41,-11.18)
\psdots[dotstyle=*](-5.28,-10.27)
\psdots[dotstyle=*](-6.52,-10.05)
\psdots[dotstyle=*](-7.66,-10.6)
\psdots[dotstyle=*](-8.25,-11.72)
\psdots[dotstyle=*](-8.07,-12.96)
\psdots[dotstyle=*](-2.31,-13.87)
\psdots[dotstyle=*](-1.07,-14.09)
\psdots[dotsize=9pt 0,dotstyle=*](0.06,-13.54)
\psdots[dotstyle=*](0.65,-12.43)
\psdots[dotstyle=*](0.48,-11.18)
\psdots[dotstyle=*](-0.39,-10.27)
\psdots[dotstyle=*](-1.63,-10.05)
\psdots[dotstyle=*](-2.77,-10.6)
\psdots[dotstyle=*](-3.36,-11.72)
\psdots[dotstyle=*](-3.18,-12.96)
\psdots[dotsize=9pt 0,dotstyle=*](2.34,-13.87)
\psdots[dotstyle=*](3.58,-14.09)
\psdots[dotstyle=*](4.71,-13.54)
\psdots[dotstyle=*](5.3,-12.43)
\psdots[dotstyle=*](5.13,-11.18)
\psdots[dotstyle=*](4.26,-10.27)
\psdots[dotstyle=*](3.02,-10.05)
\psdots[dotstyle=*](1.88,-10.6)
\psdots[dotstyle=*](1.29,-11.72)
\psdots[dotsize=9pt 0,dotstyle=*](1.47,-12.96)
\psdots[dotstyle=*](7.35,-13.87)
\psdots[dotstyle=*](8.59,-14.09)
\psdots[dotstyle=*](9.72,-13.54)
\psdots[dotsize=9pt 0,dotstyle=*](10.31,-12.43)
\psdots[dotsize=9pt 0,dotstyle=*](10.14,-11.18)
\psdots[dotstyle=*](9.27,-10.27)
\psdots[dotstyle=*](8.03,-10.05)
\psdots[dotstyle=*](6.89,-10.6)
\psdots[dotstyle=*](6.3,-11.72)
\psdots[dotstyle=*](6.48,-12.96)
\psdots[dotstyle=*](12.51,-13.87)
\psdots[dotstyle=*](13.75,-14.09)
\psdots[dotstyle=*](14.88,-13.54)
\psdots[dotstyle=*](15.47,-12.43)
\psdots[dotstyle=*](15.3,-11.18)
\psdots[dotstyle=*](14.43,-10.27)
\psdots[dotstyle=*](13.19,-10.05)
\psdots[dotsize=9pt 0,dotstyle=*](12.05,-10.6)
\psdots[dotsize=9pt 0,dotstyle=*](11.46,-11.72)
\psdots[dotstyle=*](11.64,-12.96)
\psdots[dotstyle=*](17.37,-13.87)
\psdots[dotstyle=*](18.61,-14.09)
\psdots[dotstyle=*](19.74,-13.54)
\psdots[dotstyle=*](20.33,-12.43)
\psdots[dotstyle=*](20.16,-11.18)
\psdots[dotsize=9pt 0,dotstyle=*](19.29,-10.27)
\psdots[dotsize=9pt 0,dotstyle=*](18.05,-10.05)
\psdots[dotstyle=*](16.91,-10.6)
\psdots[dotstyle=*](16.32,-11.72)
\psdots[dotstyle=*](16.5,-12.96)
\end{pspicture*}
\end{center}
\caption{Partitioning a Truncated Pentagon}
\label{polygonpartition}
\end{figure}
If $P$ is an $n$-gon with vertices swept
in the order $v_1,\ldots,v_n$, then
$Q_{v_i}$ is a line segment, $i=1,\ldots,n$;
$R_{v_1}$ and $R_{v_n}$ are empty; and $R_{v_i}$ is a point,
$i=2,\ldots,n-1$.
$Q_{v_1}\subset
H^+_{v_1}$, $Q_{v_n}\subset H^-_{v_n}$, and only the top vertex of
$Q_{v_i}$ is in $H^+_{v_i}$, $i=2,\ldots,n-1$.
In Figure~\ref{polygonpartition}, the first row shows a pentagon and
its truncation, with the sweeping to occur from bottom to top.  The
second row shows the result of Step~1, in which the pre-blocks
excluding the middle faces have been constructed. 
The third row shows the result of inserting the three middle faces
(one for each of 
$T(Q_{v_2})$,
$T(Q_{v_3})$, and
$T(Q_{v_4})$)
into the appropriate pre-blocks.
The fourth row shows the final partition---the first three pre-blocks
in row~3 are merged, because the partition of $T(Q(v_1))$, a truncated
line segment, has a 
single block consisting of one $1$-face and two $0$-faces.
The other three blocks in row~3 remain unmerged---each is induced by
the trivial partition of a single point $R_{v_i}$, $i=2,3,4$.
\item The square-based pyramid ($d=3$).
Figure~\ref{pyramidfig} shows the square-based pyramid $P$ with truncated
vertices.  The view is from above, and the vertices are swept in order
$v_1,\ldots,v_5$.
Figure~\ref{pyramidtrunc0} is the complete truncation of the
pyramid together with the
facet labels (the base octagon has label $2$).
Figure~\ref{pyramidtrunc} shows the blocks in the partition of $T(P)$.

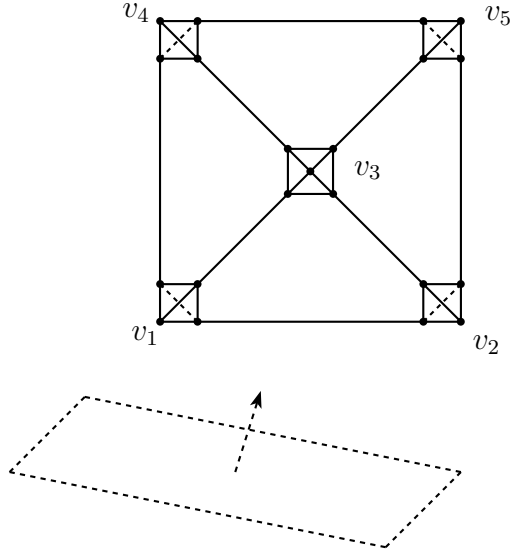
\begin{figure}
\begin{center}
\psset{xunit=1.0cm,yunit=1.0cm,algebraic=true,dotstyle=*,dotsize=3pt
0,linewidth=0.8pt,arrowsize=3pt 2,arrowinset=0.25}
\begin{pspicture*}(-4.38,-5.3)(3.7,2.9)
\psline(-2,2)(2,2)
\psline(2,2)(2,-2)
\psline(2,-2)(-2,-2)
\psline(-2,-2)(-2,2)
\psline(-2,2)(0,0)
\psline(0,0)(2,2)
\psline(2,-2)(0,0)
\psline(0,0)(-2,-2)
\psline(-2,1.5)(-1.5,1.5)
\psline(-1.5,1.5)(-1.5,2)
\psline[linestyle=dashed,dash=2pt 2pt](-2,1.5)(-1.5,2)
\psline(1.5,2)(1.5,1.5)
\psline(1.5,1.5)(2,1.5)
\psline[linestyle=dashed,dash=2pt 2pt](2,1.5)(1.5,2)
\psline(2,-1.5)(1.5,-1.5)
\psline(1.5,-1.5)(1.5,-2)
\psline[linestyle=dashed,dash=2pt 2pt](1.5,-2)(2,-1.5)
\psline(-1.5,-2)(-1.5,-1.5)
\psline(-1.5,-1.5)(-2,-1.5)
\psline[linestyle=dashed,dash=2pt 2pt](-2,-1.5)(-1.5,-2)
\psline(0.3,-0.3)(0.3,0.3)
\psline(0.3,0.3)(-0.3,0.3)
\psline(-0.3,0.3)(-0.3,-0.3)
\psline(-0.3,-0.3)(0.3,-0.3)
\psline[linestyle=dashed,dash=2pt 2pt](-3,-3)(2,-4)
\psline[linestyle=dashed,dash=2pt 2pt](-3,-3)(-4,-4)
\psline[linestyle=dashed,dash=2pt 2pt](-4,-4)(1,-5)
\psline[linestyle=dashed,dash=2pt 2pt](1,-5)(2,-4)
\psline[linestyle=dashed,dash=2pt 2pt]{->}(-1,-4)(-0.66,-2.92)
\rput[tl](-2.38,-2.06){$v_1$}
\rput[tl](2.16,-2.18){$v_2$}
\rput[tl](0.58,0.14){$v_3$}
\rput[tl](-2.5,2.24){$v_4$}
\rput[tl](2.32,2.2){$v_5$}
\psdots(0,0)
\psdots(-2,2)
\psdots(2,2)
\psdots(2,-2)
\psdots(-2,-2)
\psdots(-2,1.5)
\psdots(-1.5,1.5)
\psdots(-1.5,2)
\psdots(1.5,2)
\psdots(1.5,1.5)
\psdots(2,1.5)
\psdots(2,-1.5)
\psdots(1.5,-1.5)
\psdots(1.5,-2)
\psdots(-1.5,-2)
\psdots(-1.5,-1.5)
\psdots(-2,-1.5)
\psdots(0.3,-0.3)
\psdots(0.3,0.3)
\psdots(-0.3,0.3)
\psdots(-0.3,-0.3)
\end{pspicture*}
\end{center}
\caption{Sweeping a Pyramid (View from Above)}
\label{pyramidfig}
\end{figure}
Blocks~(1) and (2) are associated with vertex
$v_1$ of the original pyramid---note that block~(1)
also includes the truncated base of the pyramid (the outer octagon) as
well as the truncated pyramid itself.  
Block~(1) is the result of merging 9 pre-blocks, corresponding to the 9
faces in a block of the partition of $T(Q_{v_1})$ (e.g., see the first
block in the bottom row of Figure~\ref{polygonpartition}).
Block~(2) is the result of merging 4 pre-blocks, corresponding to the 4
faces in a block of the partition of $T(Q_{v_1})$ (e.g., see the second
block in the bottom row of Figure~\ref{polygonpartition}).
Neither of these
pre-blocks include middle faces, because $T(Q_{v_1})$ has none.
These two blocks are
induced by
the partition of the faces of $T(Q_{v_1})$ into two blocks.  
Blocks~(3) and (4) are
associated with vertex $v_2$.
Block~(3) is induced by the single block of the partition of
$T(Q_{v_2})$ associated with an upper vertex of $Q_{v_2}$.
Block~(4) is induced by the partition of the three faces of $T(R_{v_2})$
into a single block.  
In a similar manner, 
blocks~(5) and (6) are associated with vertex $v_3$.
Block~(7) is associated with vertex $v_4$, and is induced by the
partition of the three faces of $T(R_{v_4})$ into a single block.

\begin{figure}
\begin{center}
\psset{xunit=0.5cm,yunit=0.5cm,algebraic=true,dotstyle=*,dotsize=3pt
0,linewidth=0.8pt,arrowsize=3pt 2,arrowinset=0.25}
\begin{pspicture*}(12.24,-7.64)(23.06,2.66)
\psline(13,-5)(13,0)
\psline(13,0)(15,2)
\psline(15,2)(20,2)
\psline(20,2)(22,0)
\psline(22,0)(22,-5)
\psline(22,-5)(20,-7)
\psline(20,-7)(15,-7)
\psline(15,-7)(13,-5)
\psline(13,-5)(14,-4)
\psline(14,-4)(14,-1)
\psline(14,-1)(13,0)
\psline(14,-1)(15,-1)
\psline(15,-1)(16,0)
\psline(16,0)(16,1)
\psline(16,1)(15,2)
\psline(16,1)(19,1)
\psline(19,1)(20,2)
\psline(19,1)(19,0)
\psline(19,0)(20,-1)
\psline(20,-1)(21,-1)
\psline(21,-1)(22,0)
\psline(21,-1)(21,-4)
\psline(21,-4)(22,-5)
\psline(21,-4)(20,-4)
\psline(20,-4)(19,-5)
\psline(19,-5)(19,-6)
\psline(19,-6)(20,-7)
\psline(19,-6)(16,-6)
\psline(16,-6)(15,-7)
\psline(16,-6)(16,-5)
\psline(16,-5)(15,-4)
\psline(15,-4)(14,-4)
\psline(15,-4)(16,-3)
\psline(16,-3)(16,-2)
\psline(16,-2)(15,-1)
\psline(16,-2)(17,-1)
\psline(17,-1)(16,0)
\psline(17,-1)(18,-1)
\psline(18,-1)(19,0)
\psline(18,-1)(19,-2)
\psline(19,-2)(20,-1)
\psline(19,-2)(19,-3)
\psline(19,-3)(20,-4)
\psline(19,-5)(18,-4)
\psline(18,-4)(17,-4)
\psline(17,-4)(16,-5)
\psline(17,-4)(16,-3)
\psline(18,-4)(19,-3)
\rput[tl](17.32,1.66){$1$}
\rput[tl](15.86,-0.72){$1$}
\rput[tl](18.82,-0.78){$1$}
\rput[tl](18.8,-3.78){$1$}
\rput[tl](15.86,-3.72){$1$}
\rput[tl](14.84,-2.26){$2$}
\rput[tl](17.36,-4.78){$2$}
\rput[tl](19.82,-2.26){$2$}
\rput[tl](17.42,0.22){$2$}
\rput[tl](17.3,-2.32){$0$}
\rput[tl](14.58,-5.04){$0$}
\rput[tl](20.08,-5.08){$0$}
\rput[tl](20.16,0.56){$0$}
\rput[tl](14.52,0.5){$0$}
\rput[tl](13.36,-2.3){$1$}
\rput[tl](17.3,-6.26){$1$}
\rput[tl](21.36,-2.32){$1$}
\psdots(13,-5)
\psdots(13,0)
\psdots(15,2)
\psdots(20,2)
\psdots(22,0)
\psdots(22,-5)
\psdots(20,-7)
\psdots(15,-7)
\psdots(14,-4)
\psdots(14,-1)
\psdots(15,-1)
\psdots(16,0)
\psdots(16,1)
\psdots(19,1)
\psdots(19,0)
\psdots(20,-1)
\psdots(21,-1)
\psdots(21,-4)
\psdots(20,-4)
\psdots(19,-5)
\psdots(19,-6)
\psdots(16,-6)
\psdots(16,-5)
\psdots(15,-4)
\psdots(16,-3)
\psdots(16,-2)
\psdots(17,-1)
\psdots(18,-1)
\psdots(19,-2)
\psdots(19,-3)
\psdots(18,-4)
\psdots(17,-4)
\end{pspicture*}
\end{center}
\caption{Truncated Pyramid}
\label{pyramidtrunc0}
\end{figure}
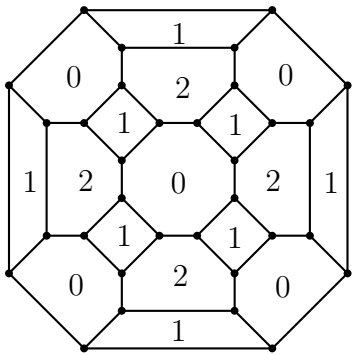

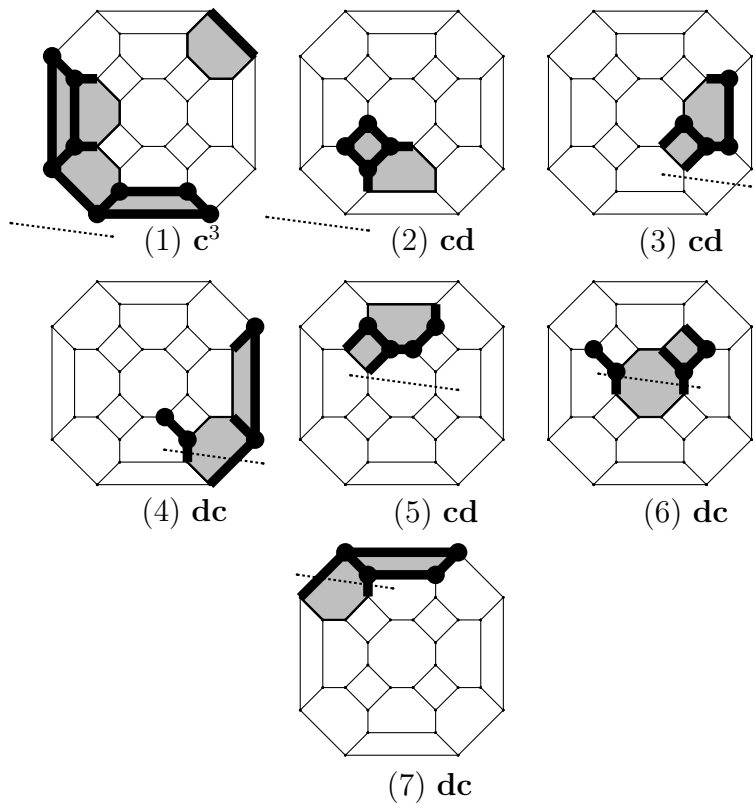
\begin{figure}
\begin{center}
\psset{xunit=0.3cm,yunit=0.3cm,algebraic=true,dotstyle=*,dotsize=3pt
0,linewidth=0.8pt,arrowsize=3pt 2,arrowinset=0.25}
\begin{pspicture*}(-0.91,-33.86)(34.79,3.95)
\pspolygon[fillcolor=lightgray,fillstyle=solid](11,0)(10,-1)(9,-1)(8,0)(8,1)(9,2)
\pspolygon[fillcolor=lightgray,fillstyle=solid](3,-4)(4,-4)(5,-3)(5,-2)(4,-1)(3,-1)
\pspolygon[fillcolor=lightgray,fillstyle=solid](3,-1)(2,0)(2,-5)(3,-4)
\pspolygon[fillcolor=lightgray,fillstyle=solid](4,-7)(5,-6)(5,-5)(4,-4)(3,-4)(2,-5)
\pspolygon[fillcolor=lightgray,fillstyle=solid](4,-7)(9,-7)(8,-6)(5,-6)
\pspolygon[fillcolor=lightgray,fillstyle=solid](16,-5)(15,-4)(16,-3)(17,-4)
\pspolygon[fillcolor=lightgray,fillstyle=solid](16,-6)(16,-5)(17,-4)(18,-4)(19,-5)(19,-6)
\pspolygon[fillcolor=lightgray,fillstyle=solid](31,-4)(30,-3)(29,-4)(30,-5)
\pspolygon[fillcolor=lightgray,fillstyle=solid](32,-4)(32,-1)(31,-1)(30,-2)(30,-3)(31,-4)
\pspolygon[fillcolor=lightgray,fillstyle=solid](11,-17)(10,-16)(9,-16)(8,-17)(8,-18)(9,-19)
\pspolygon[fillcolor=lightgray,fillstyle=solid](11,-17)(11,-12)(10,-13)(10,-16)
\pspolygon[fillcolor=lightgray,fillstyle=solid](17,-13)(16,-11.96)(15,-13)(16,-14)
\pspolygon[fillcolor=lightgray,fillstyle=solid](19,-12)(19,-11)(16,-11)(16,-11.96)(17,-13)(18,-13)
\pspolygon[fillcolor=lightgray,fillstyle=solid](30,-14)(29,-13)(28,-13)(27,-14)(27,-15)(28,-16)(29,-16)(30,-15)
\pspolygon[fillcolor=lightgray,fillstyle=solid](30,-14)(31,-13)(30,-12)(29,-13)
\pspolygon[fillcolor=lightgray,fillstyle=solid](15,-22)(13,-24)(14,-25)(15,-25)(16,-23.96)(16,-23)
\pspolygon[fillcolor=lightgray,fillstyle=solid](15,-22)(20,-22)(19,-23)(16,-23)
\psline[linewidth=0.4pt](6,-1)(7,-1)
\psline[linewidth=0.4pt](7,-1)(8,-2)
\psline[linewidth=0.4pt](8,-2)(8,-3)
\psline[linewidth=0.4pt](8,-3)(7,-4)
\psline[linewidth=0.4pt](7,-4)(6,-4)
\psline[linewidth=0.4pt](6,-4)(5,-3)
\psline[linewidth=0.4pt](5,-3)(5,-2)
\psline[linewidth=0.4pt](5,-2)(6,-1)
\psline[linewidth=0.4pt](4,-4)(5,-5)
\psline[linewidth=0.4pt](5,-5)(5,-6)
\psline[linewidth=3.6pt](4,-4)(3,-4)
\psline[linewidth=3.6pt](3,-4)(2,-5)
\psline[linewidth=3.6pt](2,-5)(4,-7)
\psline[linewidth=3.6pt](4,-7)(5,-6)
\psline[linewidth=0.4pt](4,-4)(5,-3)
\psline[linewidth=0.4pt](5,-5)(6,-4)
\psline[linewidth=0.4pt](5,-2)(4,-1)
\psline[linewidth=0.4pt](6,-1)(5,0.04)
\psline[linewidth=0.4pt](5,0.04)(4,-1)
\psline[linewidth=3.6pt](4,-1)(3,-1)
\psline[linewidth=0.4pt](5,0.04)(5,1)
\psline[linewidth=0.4pt](5,1)(4,2)
\psline[linewidth=0.4pt](4,2)(2,0)
\psline[linewidth=3.6pt](2,0)(3,-1)
\psline[linewidth=3.6pt](2,0)(2,-5)
\psline[linewidth=3.6pt](3,-4)(3,-1)
\psline[linewidth=0.4pt](8,0)(7,-1)
\psline[linewidth=0.4pt](9,-1)(8,-2)
\psline[linewidth=0.4pt](8,0)(9,-1)
\psline[linewidth=0.4pt](8,-3)(9,-4)
\psline[linewidth=0.4pt](9,-4)(8,-5)
\psline[linewidth=0.4pt](8,-5)(7,-4)
\psline[linewidth=0.4pt](8,1)(8,0)
\psline[linewidth=0.4pt](9,-1)(10,-1)
\psline[linewidth=0.4pt](9,-4)(10,-4)
\psline[linewidth=0.4pt](8,-5)(8,-6)
\psline[linewidth=0.4pt](5,1)(8,1)
\psline[linewidth=0.4pt](8,1)(9,2)
\psline[linewidth=3.6pt](9,2)(11,0)
\psline[linewidth=0.4pt](11,0)(10,-1)
\psline[linewidth=0.4pt](4,2)(9,2)
\psline[linewidth=0.4pt](10,-1)(10,-4)
\psline[linewidth=0.4pt](11,0)(11,-5)
\psline[linewidth=0.4pt](11,-5)(10,-4)
\psline[linewidth=0.4pt](11,-5)(9,-7)
\psline[linewidth=3.6pt](9,-7)(8,-6)
\psline[linewidth=3.6pt](8,-6)(5,-6)
\psline[linewidth=3.6pt](4,-7)(9,-7)
\psline[linewidth=0.4pt](17,-1)(18,-1)
\psline[linewidth=0.4pt](18,-1)(19,-2)
\psline[linewidth=0.4pt](19,-2)(19,-3)
\psline[linewidth=0.4pt](19,-3)(18,-4)
\psline[linewidth=3.6pt](18,-4)(17,-4)
\psline[linewidth=3.6pt](17,-4)(16,-3)
\psline[linewidth=0.4pt](16,-3)(16,-2)
\psline[linewidth=0.4pt](16,-2)(17,-1)
\psline[linewidth=3.6pt](15,-4)(16,-5)
\psline[linewidth=3.6pt](16,-5)(16,-6)
\psline[linewidth=0.4pt](15,-4)(14,-4)
\psline[linewidth=0.4pt](14,-4)(13,-5)
\psline[linewidth=0.4pt](13,-5)(15,-7)
\psline[linewidth=0.4pt](15,-7)(16,-6)
\psline[linewidth=3.6pt](15,-4)(16,-3)
\psline[linewidth=3.6pt](16,-5)(17,-4)
\psline[linewidth=0.4pt](16,-2)(15,-1)
\psline[linewidth=0.4pt](17,-1)(16,0.04)
\psline[linewidth=0.4pt](16,0.04)(15,-1)
\psline[linewidth=0.4pt](15,-1)(14,-1)
\psline[linewidth=0.4pt](16,0.04)(16,1)
\psline[linewidth=0.4pt](16,1)(15,2)
\psline[linewidth=0.4pt](15,2)(13,0)
\psline[linewidth=0.4pt](13,0)(14,-1)
\psline[linewidth=0.4pt](13,0)(13,-5)
\psline[linewidth=0.4pt](14,-4)(14,-1)
\psline[linewidth=0.4pt](19,0)(18,-1)
\psline[linewidth=0.4pt](20,-1)(19,-2)
\psline[linewidth=0.4pt](19,0)(20,-1)
\psline[linewidth=0.4pt](19,-3)(20,-4)
\psline[linewidth=0.4pt](20,-4)(19,-5)
\psline[linewidth=0.4pt](19,-5)(18,-4)
\psline[linewidth=0.4pt](19,1)(19,0)
\psline[linewidth=0.4pt](20,-1)(21,-1)
\psline[linewidth=0.4pt](20,-4)(21,-4)
\psline[linewidth=0.4pt](19,-5)(19,-6)
\psline[linewidth=0.4pt](16,1)(19,1)
\psline[linewidth=0.4pt](19,1)(20,2)
\psline[linewidth=0.4pt](20,2)(22,0)
\psline[linewidth=0.4pt](22,0)(21,-1)
\psline[linewidth=0.4pt](15,2)(20,2)
\psline[linewidth=0.4pt](21,-1)(21,-4)
\psline[linewidth=0.4pt](22,0)(22,-5)
\psline[linewidth=0.4pt](22,-5)(21,-4)
\psline[linewidth=0.4pt](22,-5)(20,-7)
\psline[linewidth=0.4pt](20,-7)(19,-6)
\psline[linewidth=0.4pt](19,-6)(16,-6)
\psline[linewidth=0.4pt](15,-7)(20,-7)
\psline[linewidth=0.4pt](28,-1)(29,-1)
\psline[linewidth=0.4pt](29,-1)(30,-2)
\psline[linewidth=0.4pt](30,-2)(30,-3)
\psline[linewidth=3.6pt](30,-3)(29,-4)
\psline[linewidth=0.4pt](29,-4)(28,-4)
\psline[linewidth=0.4pt](28,-4)(27,-3)
\psline[linewidth=0.4pt](27,-3)(27,-2)
\psline[linewidth=0.4pt](27,-2)(28,-1)
\psline[linewidth=0.4pt](26,-4)(27,-5)
\psline[linewidth=0.4pt](27,-5)(27,-6)
\psline[linewidth=0.4pt](26,-4)(25,-4)
\psline[linewidth=0.4pt](25,-4)(24,-5)
\psline[linewidth=0.4pt](24,-5)(26,-7)
\psline[linewidth=0.4pt](26,-7)(27,-6)
\psline[linewidth=0.4pt](26,-4)(27,-3)
\psline[linewidth=0.4pt](27,-5)(28,-4)
\psline[linewidth=0.4pt](27,-2)(26,-1)
\psline[linewidth=0.4pt](28,-1)(27,0.04)
\psline[linewidth=0.4pt](27,0.04)(26,-1)
\psline[linewidth=0.4pt](26,-1)(25,-1)
\psline[linewidth=0.4pt](27,0.04)(27,1)
\psline[linewidth=0.4pt](27,1)(26,2)
\psline[linewidth=0.4pt](26,2)(24,0)
\psline[linewidth=0.4pt](24,0)(25,-1)
\psline[linewidth=0.4pt](24,0)(24,-5)
\psline[linewidth=0.4pt](25,-4)(25,-1)
\psline[linewidth=0.4pt](30,0)(29,-1)
\psline[linewidth=0.4pt](31,-1)(30,-2)
\psline[linewidth=0.4pt](30,0)(31,-1)
\psline[linewidth=3.6pt](30,-3)(31,-4)
\psline[linewidth=3.6pt](31,-4)(30,-5)
\psline[linewidth=0.4pt](30,-5)(29,-4)
\psline[linewidth=0.4pt](30,1)(30,0)
\psline[linewidth=3.6pt](31,-1)(32,-1)
\psline[linewidth=3.6pt](31,-4)(32,-4)
\psline[linewidth=0.4pt](30,-5)(30,-6)
\psline[linewidth=0.4pt](27,1)(30,1)
\psline[linewidth=0.4pt](30,1)(31,2)
\psline[linewidth=0.4pt](31,2)(33,0)
\psline[linewidth=0.4pt](33,0)(32,-1)
\psline[linewidth=0.4pt](26,2)(31,2)
\psline[linewidth=3.6pt](32,-1)(32,-4)
\psline[linewidth=0.4pt](33,0)(33,-5)
\psline[linewidth=0.4pt](33,-5)(32,-4)
\psline[linewidth=0.4pt](33,-5)(31,-7)
\psline[linewidth=0.4pt](31,-7)(30,-6)
\psline[linewidth=0.4pt](30,-6)(27,-6)
\psline[linewidth=0.4pt](26,-7)(31,-7)
\psline[linewidth=0.4pt](6,-13)(7,-13)
\psline[linewidth=0.4pt](7,-13)(8,-14)
\psline[linewidth=0.4pt](8,-14)(8,-15)
\psline[linewidth=0.4pt](8,-15)(7,-16)
\psline[linewidth=0.4pt](7,-16)(6,-16)
\psline[linewidth=0.4pt](6,-16)(5,-15)
\psline[linewidth=0.4pt](5,-15)(5,-14)
\psline[linewidth=0.4pt](5,-14)(6,-13)
\psline[linewidth=0.4pt](4,-16)(5,-17)
\psline[linewidth=0.4pt](5,-17)(5,-18)
\psline[linewidth=0.4pt](4,-16)(3,-16)
\psline[linewidth=0.4pt](3,-16)(2,-17)
\psline[linewidth=0.4pt](2,-17)(4,-19)
\psline[linewidth=0.4pt](4,-19)(5,-18)
\psline[linewidth=0.4pt](4,-16)(5,-15)
\psline[linewidth=0.4pt](5,-17)(6,-16)
\psline[linewidth=0.4pt](5,-14)(4,-13)
\psline[linewidth=0.4pt](6,-13)(5,-11.96)
\psline[linewidth=0.4pt](5,-11.96)(4,-13)
\psline[linewidth=0.4pt](4,-13)(3,-13)
\psline[linewidth=0.4pt](5,-11.96)(5,-11)
\psline[linewidth=0.4pt](5,-11)(4,-10)
\psline[linewidth=0.4pt](4,-10)(2,-12)
\psline[linewidth=0.4pt](2,-12)(3,-13)
\psline[linewidth=0.4pt](2,-12)(2,-17)
\psline[linewidth=0.4pt](3,-16)(3,-13)
\psline[linewidth=0.4pt](8,-12)(7,-13)
\psline[linewidth=0.4pt](9,-13)(8,-14)
\psline[linewidth=0.4pt](8,-12)(9,-13)
\psline[linewidth=0.4pt](8,-15)(9,-16)
\psline[linewidth=0.4pt](9,-16)(8,-17)
\psline[linewidth=3.6pt](8,-17)(7,-16)
\psline[linewidth=0.4pt](8,-11)(8,-12)
\psline[linewidth=0.4pt](9,-13)(10,-13)
\psline[linewidth=0.4pt](9,-16)(10,-16)
\psline[linewidth=3.6pt](8,-17)(8,-18)
\psline[linewidth=0.4pt](5,-11)(8,-11)
\psline[linewidth=0.4pt](8,-11)(9,-10)
\psline[linewidth=0.4pt](9,-10)(11,-12)
\psline[linewidth=3.6pt](11,-12)(10,-13)
\psline[linewidth=0.4pt](4,-10)(9,-10)
\psline[linewidth=0.4pt](10,-13)(10,-16)
\psline[linewidth=3.6pt](11,-12)(11,-17)
\psline[linewidth=3.6pt](11,-17)(10,-16)
\psline[linewidth=3.6pt](11,-17)(9,-19)
\psline[linewidth=0.4pt](9,-19)(8,-18)
\psline[linewidth=0.4pt](8,-18)(5,-18)
\psline[linewidth=0.4pt](4,-19)(9,-19)
\psline[linewidth=3.6pt](17,-13)(18,-13)
\psline[linewidth=0.4pt](18,-13)(19,-14)
\psline[linewidth=0.4pt](19,-14)(19,-15)
\psline[linewidth=0.4pt](19,-15)(18,-16)
\psline[linewidth=0.4pt](18,-16)(17,-16)
\psline[linewidth=0.4pt](17,-16)(16,-15)
\psline[linewidth=0.4pt](16,-15)(16,-14)
\psline[linewidth=3.6pt](16,-14)(17,-13)
\psline[linewidth=0.4pt](15,-16)(16,-17)
\psline[linewidth=0.4pt](16,-17)(16,-18)
\psline[linewidth=0.4pt](15,-16)(14,-16)
\psline[linewidth=0.4pt](14,-16)(13,-17)
\psline[linewidth=0.4pt](13,-17)(15,-19)
\psline[linewidth=0.4pt](15,-19)(16,-18)
\psline[linewidth=0.4pt](15,-16)(16,-15)
\psline[linewidth=0.4pt](16,-17)(17,-16)
\psline[linewidth=0.4pt](16,-14)(15,-13)
\psline[linewidth=3.6pt](17,-13)(16,-11.96)
\psline[linewidth=3.6pt](16,-11.96)(15,-13)
\psline[linewidth=0.4pt](15,-13)(14,-13)
\psline[linewidth=0.4pt](16,-11.96)(16,-11)
\psline[linewidth=0.4pt](16,-11)(15,-10)
\psline[linewidth=0.4pt](15,-10)(13,-12)
\psline[linewidth=0.4pt](13,-12)(14,-13)
\psline[linewidth=0.4pt](13,-12)(13,-17)
\psline[linewidth=0.4pt](14,-16)(14,-13)
\psline[linewidth=3.6pt](19,-12)(18,-13)
\psline[linewidth=0.4pt](20,-13)(19,-14)
\psline[linewidth=0.4pt](19,-12)(20,-13)
\psline[linewidth=0.4pt](19,-15)(20,-16)
\psline[linewidth=0.4pt](20,-16)(19,-17)
\psline[linewidth=0.4pt](19,-17)(18,-16)
\psline[linewidth=3.6pt](19,-11)(19,-12)
\psline[linewidth=0.4pt](20,-13)(21,-13)
\psline[linewidth=0.4pt](20,-16)(21,-16)
\psline[linewidth=0.4pt](19,-17)(19,-18)
\psline[linewidth=0.4pt](16,-11)(19,-11)
\psline[linewidth=0.4pt](19,-11)(20,-10)
\psline[linewidth=0.4pt](20,-10)(22,-12)
\psline[linewidth=0.4pt](22,-12)(21,-13)
\psline[linewidth=0.4pt](15,-10)(20,-10)
\psline[linewidth=0.4pt](21,-13)(21,-16)
\psline[linewidth=0.4pt](22,-12)(22,-17)
\psline[linewidth=0.4pt](22,-17)(21,-16)
\psline[linewidth=0.4pt](22,-17)(20,-19)
\psline[linewidth=0.4pt](20,-19)(19,-18)
\psline[linewidth=0.4pt](19,-18)(16,-18)
\psline[linewidth=0.4pt](15,-19)(20,-19)
\psline[linewidth=0.4pt](28,-13)(29,-13)
\psline[linewidth=3.6pt](29,-13)(30,-14)
\psline[linewidth=3.6pt](30,-14)(30,-15)
\psline[linewidth=0.4pt](30,-15)(29,-16)
\psline[linewidth=0.4pt](29,-16)(28,-16)
\psline[linewidth=0.4pt](28,-16)(27,-15)
\psline[linewidth=3.6pt](27,-15)(27,-14)
\psline[linewidth=0.4pt](27,-14)(28,-13)
\psline[linewidth=0.4pt](26,-16)(27,-17)
\psline[linewidth=0.4pt](27,-17)(27,-18)
\psline[linewidth=0.4pt](26,-16)(25,-16)
\psline[linewidth=0.4pt](25,-16)(24,-17)
\psline[linewidth=0.4pt](24,-17)(26,-19)
\psline[linewidth=0.4pt](26,-19)(27,-18)
\psline[linewidth=0.4pt](26,-16)(27,-15)
\psline[linewidth=0.4pt](27,-17)(28,-16)
\psline[linewidth=3.6pt](27,-14)(26,-13)
\psline[linewidth=0.4pt](28,-13)(27,-11.96)
\psline[linewidth=0.4pt](27,-11.96)(26,-13)
\psline[linewidth=0.4pt](26,-13)(25,-13)
\psline[linewidth=0.4pt](27,-11.96)(27,-11)
\psline[linewidth=0.4pt](27,-11)(26,-10)
\psline[linewidth=0.4pt](26,-10)(24,-12)
\psline[linewidth=0.4pt](24,-12)(25,-13)
\psline[linewidth=0.4pt](24,-12)(24,-17)
\psline[linewidth=0.4pt](25,-16)(25,-13)
\psline[linewidth=0.4pt](30,-12)(29,-13)
\psline[linewidth=3.6pt](31,-13)(30,-14)
\psline[linewidth=3.6pt](30,-12)(31,-13)
\psline[linewidth=0.4pt](30,-15)(31,-16)
\psline[linewidth=0.4pt](31,-16)(30,-17)
\psline[linewidth=0.4pt](30,-17)(29,-16)
\psline[linewidth=0.4pt](30,-11)(30,-12)
\psline[linewidth=0.4pt](31,-13)(32,-13)
\psline[linewidth=0.4pt](31,-16)(32,-16)
\psline[linewidth=0.4pt](30,-17)(30,-18)
\psline[linewidth=0.4pt](27,-11)(30,-11)
\psline[linewidth=0.4pt](30,-11)(31,-10)
\psline[linewidth=0.4pt](31,-10)(33,-12)
\psline[linewidth=0.4pt](33,-12)(32,-13)
\psline[linewidth=0.4pt](26,-10)(31,-10)
\psline[linewidth=0.4pt](32,-13)(32,-16)
\psline[linewidth=0.4pt](33,-12)(33,-17)
\psline[linewidth=0.4pt](33,-17)(32,-16)
\psline[linewidth=0.4pt](33,-17)(31,-19)
\psline[linewidth=0.4pt](31,-19)(30,-18)
\psline[linewidth=0.4pt](30,-18)(27,-18)
\psline[linewidth=0.4pt](26,-19)(31,-19)
\psline[linewidth=0.4pt](17,-25)(18,-25)
\psline[linewidth=0.4pt](18,-25)(19,-26)
\psline[linewidth=0.4pt](19,-26)(19,-27)
\psline[linewidth=0.4pt](19,-27)(18,-28)
\psline[linewidth=0.4pt](18,-28)(17,-28)
\psline[linewidth=0.4pt](17,-28)(16,-27)
\psline[linewidth=0.4pt](16,-27)(16,-26)
\psline[linewidth=0.4pt](16,-26)(17,-25)
\psline[linewidth=0.4pt](15,-28)(16,-29)
\psline[linewidth=0.4pt](16,-29)(16,-30)
\psline[linewidth=0.4pt](15,-28)(14,-28)
\psline[linewidth=0.4pt](14,-28)(13,-29)
\psline[linewidth=0.4pt](13,-29)(15,-31)
\psline[linewidth=0.4pt](15,-31)(16,-30)
\psline[linewidth=0.4pt](15,-28)(16,-27)
\psline[linewidth=0.4pt](16,-29)(17,-28)
\psline[linewidth=0.4pt](16,-26)(15,-25)
\psline[linewidth=0.4pt](17,-25)(16,-23.96)
\psline[linewidth=0.4pt](16,-23.96)(15,-25)
\psline[linewidth=0.4pt](15,-25)(14,-25)
\psline[linewidth=3.6pt](16,-23.96)(16,-23)
\psline[linewidth=3.6pt](16,-23)(15,-22)
\psline[linewidth=3.6pt](15,-22)(13,-24)
\psline[linewidth=0.4pt](13,-24)(14,-25)
\psline[linewidth=0.4pt](13,-24)(13,-29)
\psline[linewidth=0.4pt](14,-28)(14,-25)
\psline[linewidth=0.4pt](19,-24)(18,-25)
\psline[linewidth=0.4pt](20,-25)(19,-26)
\psline[linewidth=0.4pt](19,-24)(20,-25)
\psline[linewidth=0.4pt](19,-27)(20,-28)
\psline[linewidth=0.4pt](20,-28)(19,-29)
\psline[linewidth=0.4pt](19,-29)(18,-28)
\psline[linewidth=0.4pt](19,-23)(19,-24)
\psline[linewidth=0.4pt](20,-25)(21,-25)
\psline[linewidth=0.4pt](20,-28)(21,-28)
\psline[linewidth=0.4pt](19,-29)(19,-30)
\psline[linewidth=3.6pt](16,-23)(19,-23)
\psline[linewidth=3.6pt](19,-23)(20,-22)
\psline[linewidth=0.4pt](20,-22)(22,-24)
\psline[linewidth=0.4pt](22,-24)(21,-25)
\psline[linewidth=3.6pt](15,-22)(20,-22)
\psline[linewidth=0.4pt](21,-25)(21,-28)
\psline[linewidth=0.4pt](22,-24)(22,-29)
\psline[linewidth=0.4pt](22,-29)(21,-28)
\psline[linewidth=0.4pt](22,-29)(20,-31)
\psline[linewidth=0.4pt](20,-31)(19,-30)
\psline[linewidth=0.4pt](19,-30)(16,-30)
\psline[linewidth=0.4pt](15,-31)(20,-31)
\psline(11,0)(10,-1)
\psline(10,-1)(9,-1)
\psline(9,-1)(8,0)
\psline(8,0)(8,1)
\psline(8,1)(9,2)
\psline(9,2)(11,0)
\psline(3,-4)(4,-4)
\psline(4,-4)(5,-3)
\psline(5,-3)(5,-2)
\psline(5,-2)(4,-1)
\psline(4,-1)(3,-1)
\psline(3,-1)(3,-4)
\psline(3,-1)(2,0)
\psline(2,0)(2,-5)
\psline(2,-5)(3,-4)
\psline(3,-4)(3,-1)
\psline(4,-7)(5,-6)
\psline(5,-6)(5,-5)
\psline(5,-5)(4,-4)
\psline(4,-4)(3,-4)
\psline(3,-4)(2,-5)
\psline(2,-5)(4,-7)
\psline(4,-7)(9,-7)
\psline(9,-7)(8,-6)
\psline(8,-6)(5,-6)
\psline(5,-6)(4,-7)
\psline[linestyle=dashed,dash=1pt 1pt](11.51,-7.13)(16.03,-7.73)
\psline[linestyle=dashed,dash=1pt 1pt](29.06,-5.24)(32.97,-5.78)
\psline[linestyle=dashed,dash=1pt 1pt](6.98,-17.45)(11.37,-18.06)
\psline[linestyle=dashed,dash=1pt 1pt](26.22,-14.08)(30.68,-14.69)
\psline[linestyle=dashed,dash=1pt 1pt](0.17,-7.4)(4.69,-8)
\psline[linestyle=dashed,dash=1pt 1pt](15.22,-14.15)(20,-14.8)
\psline[linestyle=dashed,dash=1pt 1pt](12.86,-22.99)(17.11,-23.6)
\rput[tl](6.04,-7.50){$(1)\ \bc^3$}
\psline(16,-5)(15,-4)
\psline(15,-4)(16,-3)
\psline(16,-3)(17,-4)
\psline(17,-4)(16,-5)
\psline(16,-6)(16,-5)
\psline(16,-5)(17,-4)
\psline(17,-4)(18,-4)
\psline(18,-4)(19,-5)
\psline(19,-5)(19,-6)
\psline(19,-6)(16,-6)
\rput[tl](16.97,-7.53){$(2)\ \bc\bd$}
\psline(31,-4)(30,-3)
\psline(30,-3)(29,-4)
\psline(29,-4)(30,-5)
\psline(30,-5)(31,-4)
\psline(32,-4)(32,-1)
\psline(32,-1)(31,-1)
\psline(31,-1)(30,-2)
\psline(30,-2)(30,-3)
\psline(30,-3)(31,-4)
\psline(31,-4)(32,-4)
\rput[tl](27.98,-7.6){$(3)\ \bc\bd$}
\psline(11,-17)(10,-16)
\psline(10,-16)(9,-16)
\psline(9,-16)(8,-17)
\psline(8,-17)(8,-18)
\psline(8,-18)(9,-19)
\psline(9,-19)(11,-17)
\psline(11,-17)(11,-12)
\psline(11,-12)(10,-13)
\psline(10,-13)(10,-16)
\psline(10,-16)(11,-17)
\rput[tl](5.97,-19.55){$(4)\ \bd\bc$}
\psline(17,-13)(16,-11.96)
\psline(16,-11.96)(15,-13)
\psline(15,-13)(16,-14)
\psline(16,-14)(17,-13)
\psline(19,-12)(19,-11)
\psline(19,-11)(16,-11)
\psline(16,-11)(16,-11.96)
\psline(16,-11.96)(17,-13)
\psline(17,-13)(18,-13)
\psline(18,-13)(19,-12)
\rput[tl](17.11,-19.61){$(5)\ \bcd$}
\psline(30,-14)(29,-13)
\psline(29,-13)(28,-13)
\psline(28,-13)(27,-14)
\psline(27,-14)(27,-15)
\psline(27,-15)(28,-16)
\psline(28,-16)(29,-16)
\psline(29,-16)(30,-15)
\psline(30,-15)(30,-14)
\psline(30,-14)(31,-13)
\psline(31,-13)(30,-12)
\psline(30,-12)(29,-13)
\psline(29,-13)(30,-14)
\rput[tl](28.18,-19.61){$(6)\ \bd\bc$}
\psline(15,-22)(13,-24)
\psline(13,-24)(14,-25)
\psline(14,-25)(15,-25)
\psline(15,-25)(16,-23.96)
\psline(16,-23.96)(16,-23)
\psline(16,-23)(15,-22)
\psline(15,-22)(20,-22)
\psline(20,-22)(19,-23)
\psline(19,-23)(16,-23)
\psline(16,-23)(15,-22)
\rput[tl](16.84,-31.7){$(7)\ \bd\bc$}
\psdots[dotsize=1pt 0](6,-1)
\psdots[dotsize=1pt 0](7,-1)
\psdots[dotsize=1pt 0](8,-2)
\psdots[dotsize=1pt 0](8,-3)
\psdots[dotsize=1pt 0](7,-4)
\psdots[dotsize=1pt 0](6,-4)
\psdots[dotsize=1pt 0](5,-3)
\psdots[dotsize=1pt 0](5,-2)
\psdots[dotsize=1pt 0](4,-4)
\psdots[dotsize=1pt 0](5,-5)
\psdots[dotsize=7pt 0](5,-6)
\psdots[dotsize=7pt 0](3,-4)
\psdots[dotsize=7pt 0](2,-5)
\psdots[dotsize=7pt 0](4,-7)
\psdots[dotsize=1pt 0](4,-1)
\psdots[dotsize=1pt 0](5,0.04)
\psdots[dotsize=7pt 0](3,-1)
\psdots[dotsize=1pt 0](5,1)
\psdots[dotsize=1pt 0](4,2)
\psdots[dotsize=7pt 0](2,0)
\psdots[dotsize=1pt 0](8,0)
\psdots[dotsize=1pt 0](9,-1)
\psdots[dotsize=1pt 0](9,-4)
\psdots[dotsize=1pt 0](8,-5)
\psdots[dotsize=1pt 0](8,1)
\psdots[dotsize=1pt 0](10,-1)
\psdots[dotsize=1pt 0](10,-4)
\psdots[dotsize=7pt 0](8,-6)
\psdots[dotsize=1pt 0](9,2)
\psdots[dotsize=1pt 0](11,0)
\psdots[dotsize=1pt 0](11,-5)
\psdots[dotsize=7pt 0](9,-7)
\psdots[dotsize=1pt 0](17,-1)
\psdots[dotsize=1pt 0](18,-1)
\psdots[dotsize=1pt 0](18,-1)
\psdots[dotsize=1pt 0](19,-2)
\psdots[dotsize=1pt 0](19,-2)
\psdots[dotsize=1pt 0](19,-3)
\psdots[dotsize=1pt 0](19,-3)
\psdots[dotsize=1pt 0](18,-4)
\psdots[dotsize=1pt 0](18,-4)
\psdots[dotsize=7pt 0](17,-4)
\psdots[dotsize=1pt 0](17,-4)
\psdots[dotsize=7pt 0](16,-3)
\psdots[dotsize=1pt 0](16,-3)
\psdots[dotsize=1pt 0](16,-2)
\psdots[dotsize=1pt 0](16,-2)
\psdots[dotsize=1pt 0](17,-1)
\psdots[dotsize=7pt 0](15,-4)
\psdots[dotsize=7pt 0](16,-5)
\psdots[dotsize=1pt 0](16,-5)
\psdots[dotsize=1pt 0](16,-6)
\psdots[dotsize=1pt 0](15,-4)
\psdots[dotsize=1pt 0](14,-4)
\psdots[dotsize=1pt 0](14,-4)
\psdots[dotsize=1pt 0](13,-5)
\psdots[dotsize=1pt 0](13,-5)
\psdots[dotsize=1pt 0](15,-7)
\psdots[dotsize=1pt 0](15,-7)
\psdots[dotsize=1pt 0](16,-6)
\psdots[dotsize=1pt 0](15,-4)
\psdots[dotsize=1pt 0](16,-3)
\psdots[dotsize=1pt 0](16,-5)
\psdots[dotsize=1pt 0](17,-4)
\psdots[dotsize=1pt 0](16,-2)
\psdots[dotsize=1pt 0](15,-1)
\psdots[dotsize=1pt 0](17,-1)
\psdots[dotsize=1pt 0](16,0.04)
\psdots[dotsize=1pt 0](16,0.04)
\psdots[dotsize=1pt 0](15,-1)
\psdots[dotsize=1pt 0](15,-1)
\psdots[dotsize=1pt 0](14,-1)
\psdots[dotsize=1pt 0](16,0.04)
\psdots[dotsize=1pt 0](16,1)
\psdots[dotsize=1pt 0](16,1)
\psdots[dotsize=1pt 0](15,2)
\psdots[dotsize=1pt 0](15,2)
\psdots[dotsize=1pt 0](13,0)
\psdots[dotsize=1pt 0](13,0)
\psdots[dotsize=1pt 0](14,-1)
\psdots[dotsize=1pt 0](13,0)
\psdots[dotsize=1pt 0](13,-5)
\psdots[dotsize=1pt 0](14,-4)
\psdots[dotsize=1pt 0](14,-1)
\psdots[dotsize=1pt 0](19,0)
\psdots[dotsize=1pt 0](18,-1)
\psdots[dotsize=1pt 0](20,-1)
\psdots[dotsize=1pt 0](19,-2)
\psdots[dotsize=1pt 0](19,0)
\psdots[dotsize=1pt 0](20,-1)
\psdots[dotsize=1pt 0](19,-3)
\psdots[dotsize=1pt 0](20,-4)
\psdots[dotsize=1pt 0](20,-4)
\psdots[dotsize=1pt 0](19,-5)
\psdots[dotsize=1pt 0](19,-5)
\psdots[dotsize=1pt 0](18,-4)
\psdots[dotsize=1pt 0](19,1)
\psdots[dotsize=1pt 0](19,0)
\psdots[dotsize=1pt 0](20,-1)
\psdots[dotsize=1pt 0](21,-1)
\psdots[dotsize=1pt 0](20,-4)
\psdots[dotsize=1pt 0](21,-4)
\psdots[dotsize=1pt 0](19,-5)
\psdots[dotsize=1pt 0](19,-6)
\psdots[dotsize=1pt 0](16,1)
\psdots[dotsize=1pt 0](19,1)
\psdots[dotsize=1pt 0](19,1)
\psdots[dotsize=1pt 0](20,2)
\psdots[dotsize=1pt 0](20,2)
\psdots[dotsize=1pt 0](22,0)
\psdots[dotsize=1pt 0](22,0)
\psdots[dotsize=1pt 0](21,-1)
\psdots[dotsize=1pt 0](15,2)
\psdots[dotsize=1pt 0](20,2)
\psdots[dotsize=1pt 0](21,-1)
\psdots[dotsize=1pt 0](21,-4)
\psdots[dotsize=1pt 0](22,0)
\psdots[dotsize=1pt 0](22,-5)
\psdots[dotsize=1pt 0](22,-5)
\psdots[dotsize=1pt 0](21,-4)
\psdots[dotsize=1pt 0](22,-5)
\psdots[dotsize=1pt 0](20,-7)
\psdots[dotsize=1pt 0](20,-7)
\psdots[dotsize=1pt 0](19,-6)
\psdots[dotsize=1pt 0](19,-6)
\psdots[dotsize=1pt 0](16,-6)
\psdots[dotsize=1pt 0](15,-7)
\psdots[dotsize=1pt 0](20,-7)
\psdots[dotsize=1pt 0](28,-1)
\psdots[dotsize=1pt 0](29,-1)
\psdots[dotsize=1pt 0](29,-1)
\psdots[dotsize=1pt 0](30,-2)
\psdots[dotsize=1pt 0](30,-2)
\psdots[dotsize=7pt 0](30,-3)
\psdots[dotsize=1pt 0](30,-3)
\psdots[dotsize=1pt 0](29,-4)
\psdots[dotsize=1pt 0](29,-4)
\psdots[dotsize=1pt 0](28,-4)
\psdots[dotsize=1pt 0](28,-4)
\psdots[dotsize=1pt 0](27,-3)
\psdots[dotsize=1pt 0](27,-3)
\psdots[dotsize=1pt 0](27,-2)
\psdots[dotsize=1pt 0](27,-2)
\psdots[dotsize=1pt 0](28,-1)
\psdots[dotsize=1pt 0](26,-4)
\psdots[dotsize=1pt 0](27,-5)
\psdots[dotsize=1pt 0](27,-5)
\psdots[dotsize=1pt 0](27,-6)
\psdots[dotsize=1pt 0](26,-4)
\psdots[dotsize=1pt 0](25,-4)
\psdots[dotsize=1pt 0](25,-4)
\psdots[dotsize=1pt 0](24,-5)
\psdots[dotsize=1pt 0](24,-5)
\psdots[dotsize=1pt 0](26,-7)
\psdots[dotsize=1pt 0](26,-7)
\psdots[dotsize=1pt 0](27,-6)
\psdots[dotsize=1pt 0](26,-4)
\psdots[dotsize=1pt 0](27,-3)
\psdots[dotsize=1pt 0](27,-5)
\psdots[dotsize=1pt 0](28,-4)
\psdots[dotsize=1pt 0](27,-2)
\psdots[dotsize=1pt 0](26,-1)
\psdots[dotsize=1pt 0](28,-1)
\psdots[dotsize=1pt 0](27,0.04)
\psdots[dotsize=1pt 0](27,0.04)
\psdots[dotsize=1pt 0](26,-1)
\psdots[dotsize=1pt 0](26,-1)
\psdots[dotsize=1pt 0](25,-1)
\psdots[dotsize=1pt 0](27,0.04)
\psdots[dotsize=1pt 0](27,1)
\psdots[dotsize=1pt 0](27,1)
\psdots[dotsize=1pt 0](26,2)
\psdots[dotsize=1pt 0](26,2)
\psdots[dotsize=1pt 0](24,0)
\psdots[dotsize=1pt 0](24,0)
\psdots[dotsize=1pt 0](25,-1)
\psdots[dotsize=1pt 0](24,0)
\psdots[dotsize=1pt 0](24,-5)
\psdots[dotsize=1pt 0](25,-4)
\psdots[dotsize=1pt 0](25,-1)
\psdots[dotsize=1pt 0](30,0)
\psdots[dotsize=1pt 0](29,-1)
\psdots[dotsize=1pt 0](31,-1)
\psdots[dotsize=1pt 0](30,-2)
\psdots[dotsize=1pt 0](30,0)
\psdots[dotsize=1pt 0](31,-1)
\psdots[dotsize=1pt 0](30,-3)
\psdots[dotsize=7pt 0](31,-4)
\psdots[dotsize=1pt 0](31,-4)
\psdots[dotsize=1pt 0](30,-5)
\psdots[dotsize=1pt 0](30,-5)
\psdots[dotsize=1pt 0](29,-4)
\psdots[dotsize=1pt 0](30,1)
\psdots[dotsize=1pt 0](30,0)
\psdots[dotsize=1pt 0](31,-1)
\psdots[dotsize=7pt 0](32,-1)
\psdots[dotsize=1pt 0](31,-4)
\psdots[dotsize=7pt 0](32,-4)
\psdots[dotsize=1pt 0](30,-5)
\psdots[dotsize=1pt 0](30,-6)
\psdots[dotsize=1pt 0](27,1)
\psdots[dotsize=1pt 0](30,1)
\psdots[dotsize=1pt 0](30,1)
\psdots[dotsize=1pt 0](31,2)
\psdots[dotsize=1pt 0](31,2)
\psdots[dotsize=1pt 0](33,0)
\psdots[dotsize=1pt 0](33,0)
\psdots[dotsize=1pt 0](32,-1)
\psdots[dotsize=1pt 0](26,2)
\psdots[dotsize=1pt 0](31,2)
\psdots[dotsize=1pt 0](32,-1)
\psdots[dotsize=1pt 0](32,-4)
\psdots[dotsize=1pt 0](33,0)
\psdots[dotsize=1pt 0](33,-5)
\psdots[dotsize=1pt 0](33,-5)
\psdots[dotsize=1pt 0](32,-4)
\psdots[dotsize=1pt 0](33,-5)
\psdots[dotsize=1pt 0](31,-7)
\psdots[dotsize=1pt 0](31,-7)
\psdots[dotsize=1pt 0](30,-6)
\psdots[dotsize=1pt 0](30,-6)
\psdots[dotsize=1pt 0](27,-6)
\psdots[dotsize=1pt 0](26,-7)
\psdots[dotsize=1pt 0](31,-7)
\psdots[dotsize=1pt 0](6,-13)
\psdots[dotsize=1pt 0](7,-13)
\psdots[dotsize=1pt 0](7,-13)
\psdots[dotsize=1pt 0](8,-14)
\psdots[dotsize=1pt 0](8,-14)
\psdots[dotsize=1pt 0](8,-15)
\psdots[dotsize=1pt 0](8,-15)
\psdots[dotsize=7pt 0](7,-16)
\psdots[dotsize=1pt 0](7,-16)
\psdots[dotsize=1pt 0](6,-16)
\psdots[dotsize=1pt 0](6,-16)
\psdots[dotsize=1pt 0](5,-15)
\psdots[dotsize=1pt 0](5,-15)
\psdots[dotsize=1pt 0](5,-14)
\psdots[dotsize=1pt 0](5,-14)
\psdots[dotsize=1pt 0](6,-13)
\psdots[dotsize=1pt 0](4,-16)
\psdots[dotsize=1pt 0](5,-17)
\psdots[dotsize=1pt 0](5,-17)
\psdots[dotsize=1pt 0](5,-18)
\psdots[dotsize=1pt 0](4,-16)
\psdots[dotsize=1pt 0](3,-16)
\psdots[dotsize=1pt 0](3,-16)
\psdots[dotsize=1pt 0](2,-17)
\psdots[dotsize=1pt 0](2,-17)
\psdots[dotsize=1pt 0](4,-19)
\psdots[dotsize=1pt 0](4,-19)
\psdots[dotsize=1pt 0](5,-18)
\psdots[dotsize=1pt 0](4,-16)
\psdots[dotsize=1pt 0](5,-15)
\psdots[dotsize=1pt 0](5,-17)
\psdots[dotsize=1pt 0](6,-16)
\psdots[dotsize=1pt 0](5,-14)
\psdots[dotsize=1pt 0](4,-13)
\psdots[dotsize=1pt 0](6,-13)
\psdots[dotsize=1pt 0](5,-11.96)
\psdots[dotsize=1pt 0](5,-11.96)
\psdots[dotsize=1pt 0](4,-13)
\psdots[dotsize=1pt 0](4,-13)
\psdots[dotsize=1pt 0](3,-13)
\psdots[dotsize=1pt 0](5,-11.96)
\psdots[dotsize=1pt 0](5,-11)
\psdots[dotsize=1pt 0](5,-11)
\psdots[dotsize=1pt 0](4,-10)
\psdots[dotsize=1pt 0](4,-10)
\psdots[dotsize=1pt 0](2,-12)
\psdots[dotsize=1pt 0](2,-12)
\psdots[dotsize=1pt 0](3,-13)
\psdots[dotsize=1pt 0](2,-12)
\psdots[dotsize=1pt 0](2,-17)
\psdots[dotsize=1pt 0](3,-16)
\psdots[dotsize=1pt 0](3,-13)
\psdots[dotsize=1pt 0](8,-12)
\psdots[dotsize=1pt 0](7,-13)
\psdots[dotsize=1pt 0](9,-13)
\psdots[dotsize=1pt 0](8,-14)
\psdots[dotsize=1pt 0](8,-12)
\psdots[dotsize=1pt 0](9,-13)
\psdots[dotsize=1pt 0](8,-15)
\psdots[dotsize=1pt 0](9,-16)
\psdots[dotsize=1pt 0](9,-16)
\psdots[dotsize=7pt 0](8,-17)
\psdots[dotsize=1pt 0](8,-17)
\psdots[dotsize=1pt 0](7,-16)
\psdots[dotsize=1pt 0](8,-11)
\psdots[dotsize=1pt 0](8,-12)
\psdots[dotsize=1pt 0](9,-13)
\psdots[dotsize=1pt 0](10,-13)
\psdots[dotsize=1pt 0](9,-16)
\psdots[dotsize=1pt 0](10,-16)
\psdots[dotsize=1pt 0](8,-17)
\psdots[dotsize=1pt 0](8,-18)
\psdots[dotsize=1pt 0](5,-11)
\psdots[dotsize=1pt 0](8,-11)
\psdots[dotsize=1pt 0](8,-11)
\psdots[dotsize=1pt 0](9,-10)
\psdots[dotsize=1pt 0](9,-10)
\psdots[dotsize=7pt 0](11,-12)
\psdots[dotsize=1pt 0](11,-12)
\psdots[dotsize=1pt 0](10,-13)
\psdots[dotsize=1pt 0](4,-10)
\psdots[dotsize=1pt 0](9,-10)
\psdots[dotsize=1pt 0](10,-13)
\psdots[dotsize=1pt 0](10,-16)
\psdots[dotsize=1pt 0](11,-12)
\psdots[dotsize=7pt 0](11,-17)
\psdots[dotsize=1pt 0](11,-17)
\psdots[dotsize=1pt 0](10,-16)
\psdots[dotsize=1pt 0](11,-17)
\psdots[dotsize=1pt 0](9,-19)
\psdots[dotsize=1pt 0](9,-19)
\psdots[dotsize=1pt 0](8,-18)
\psdots[dotsize=1pt 0](8,-18)
\psdots[dotsize=1pt 0](5,-18)
\psdots[dotsize=1pt 0](4,-19)
\psdots[dotsize=1pt 0](9,-19)
\psdots[dotsize=7pt 0](17,-13)
\psdots[dotsize=7pt 0](18,-13)
\psdots[dotsize=1pt 0](18,-13)
\psdots[dotsize=1pt 0](19,-14)
\psdots[dotsize=1pt 0](19,-14)
\psdots[dotsize=1pt 0](19,-15)
\psdots[dotsize=1pt 0](19,-15)
\psdots[dotsize=1pt 0](18,-16)
\psdots[dotsize=1pt 0](18,-16)
\psdots[dotsize=1pt 0](17,-16)
\psdots[dotsize=1pt 0](17,-16)
\psdots[dotsize=1pt 0](16,-15)
\psdots[dotsize=1pt 0](16,-15)
\psdots[dotsize=1pt 0](16,-14)
\psdots[dotsize=1pt 0](16,-14)
\psdots[dotsize=1pt 0](17,-13)
\psdots[dotsize=1pt 0](15,-16)
\psdots[dotsize=1pt 0](16,-17)
\psdots[dotsize=1pt 0](16,-17)
\psdots[dotsize=1pt 0](16,-18)
\psdots[dotsize=1pt 0](15,-16)
\psdots[dotsize=1pt 0](14,-16)
\psdots[dotsize=1pt 0](14,-16)
\psdots[dotsize=1pt 0](13,-17)
\psdots[dotsize=1pt 0](13,-17)
\psdots[dotsize=1pt 0](15,-19)
\psdots[dotsize=1pt 0](15,-19)
\psdots[dotsize=1pt 0](16,-18)
\psdots[dotsize=1pt 0](15,-16)
\psdots[dotsize=1pt 0](16,-15)
\psdots[dotsize=1pt 0](16,-17)
\psdots[dotsize=1pt 0](17,-16)
\psdots[dotsize=1pt 0](16,-14)
\psdots[dotsize=1pt 0](15,-13)
\psdots[dotsize=1pt 0](17,-13)
\psdots[dotsize=7pt 0](16,-11.96)
\psdots[dotsize=1pt 0](16,-11.96)
\psdots[dotsize=1pt 0](15,-13)
\psdots[dotsize=1pt 0](15,-13)
\psdots[dotsize=1pt 0](14,-13)
\psdots[dotsize=1pt 0](16,-11.96)
\psdots[dotsize=1pt 0](16,-11)
\psdots[dotsize=1pt 0](16,-11)
\psdots[dotsize=1pt 0](15,-10)
\psdots[dotsize=1pt 0](15,-10)
\psdots[dotsize=1pt 0](13,-12)
\psdots[dotsize=1pt 0](13,-12)
\psdots[dotsize=1pt 0](14,-13)
\psdots[dotsize=1pt 0](13,-12)
\psdots[dotsize=1pt 0](13,-17)
\psdots[dotsize=1pt 0](14,-16)
\psdots[dotsize=1pt 0](14,-13)
\psdots[dotsize=7pt 0](19,-12)
\psdots[dotsize=1pt 0](18,-13)
\psdots[dotsize=1pt 0](20,-13)
\psdots[dotsize=1pt 0](19,-14)
\psdots[dotsize=1pt 0](19,-12)
\psdots[dotsize=1pt 0](20,-13)
\psdots[dotsize=1pt 0](19,-15)
\psdots[dotsize=1pt 0](20,-16)
\psdots[dotsize=1pt 0](20,-16)
\psdots[dotsize=1pt 0](19,-17)
\psdots[dotsize=1pt 0](19,-17)
\psdots[dotsize=1pt 0](18,-16)
\psdots[dotsize=1pt 0](19,-11)
\psdots[dotsize=1pt 0](19,-12)
\psdots[dotsize=1pt 0](20,-13)
\psdots[dotsize=1pt 0](21,-13)
\psdots[dotsize=1pt 0](20,-16)
\psdots[dotsize=1pt 0](21,-16)
\psdots[dotsize=1pt 0](19,-17)
\psdots[dotsize=1pt 0](19,-18)
\psdots[dotsize=1pt 0](16,-11)
\psdots[dotsize=1pt 0](19,-11)
\psdots[dotsize=1pt 0](19,-11)
\psdots[dotsize=1pt 0](20,-10)
\psdots[dotsize=1pt 0](20,-10)
\psdots[dotsize=1pt 0](22,-12)
\psdots[dotsize=1pt 0](22,-12)
\psdots[dotsize=1pt 0](21,-13)
\psdots[dotsize=1pt 0](15,-10)
\psdots[dotsize=1pt 0](20,-10)
\psdots[dotsize=1pt 0](21,-13)
\psdots[dotsize=1pt 0](21,-16)
\psdots[dotsize=1pt 0](22,-12)
\psdots[dotsize=1pt 0](22,-17)
\psdots[dotsize=1pt 0](22,-17)
\psdots[dotsize=1pt 0](21,-16)
\psdots[dotsize=1pt 0](22,-17)
\psdots[dotsize=1pt 0](20,-19)
\psdots[dotsize=1pt 0](20,-19)
\psdots[dotsize=1pt 0](19,-18)
\psdots[dotsize=1pt 0](19,-18)
\psdots[dotsize=1pt 0](16,-18)
\psdots[dotsize=1pt 0](15,-19)
\psdots[dotsize=1pt 0](20,-19)
\psdots[dotsize=1pt 0](28,-13)
\psdots[dotsize=1pt 0](29,-13)
\psdots[dotsize=1pt 0](29,-13)
\psdots[dotsize=7pt 0](30,-14)
\psdots[dotsize=1pt 0](30,-14)
\psdots[dotsize=1pt 0](30,-15)
\psdots[dotsize=1pt 0](30,-15)
\psdots[dotsize=1pt 0](29,-16)
\psdots[dotsize=1pt 0](29,-16)
\psdots[dotsize=1pt 0](28,-16)
\psdots[dotsize=1pt 0](28,-16)
\psdots[dotsize=1pt 0](27,-15)
\psdots[dotsize=1pt 0](27,-15)
\psdots[dotsize=7pt 0](27,-14)
\psdots[dotsize=1pt 0](27,-14)
\psdots[dotsize=1pt 0](28,-13)
\psdots[dotsize=1pt 0](26,-16)
\psdots[dotsize=1pt 0](27,-17)
\psdots[dotsize=1pt 0](27,-17)
\psdots[dotsize=1pt 0](27,-18)
\psdots[dotsize=1pt 0](26,-16)
\psdots[dotsize=1pt 0](25,-16)
\psdots[dotsize=1pt 0](25,-16)
\psdots[dotsize=1pt 0](24,-17)
\psdots[dotsize=1pt 0](24,-17)
\psdots[dotsize=1pt 0](26,-19)
\psdots[dotsize=1pt 0](26,-19)
\psdots[dotsize=1pt 0](27,-18)
\psdots[dotsize=1pt 0](26,-16)
\psdots[dotsize=1pt 0](27,-15)
\psdots[dotsize=1pt 0](27,-17)
\psdots[dotsize=1pt 0](28,-16)
\psdots[dotsize=1pt 0](27,-14)
\psdots[dotsize=7pt 0](26,-13)
\psdots[dotsize=1pt 0](28,-13)
\psdots[dotsize=1pt 0](27,-11.96)
\psdots[dotsize=1pt 0](27,-11.96)
\psdots[dotsize=1pt 0](26,-13)
\psdots[dotsize=1pt 0](26,-13)
\psdots[dotsize=1pt 0](25,-13)
\psdots[dotsize=1pt 0](27,-11.96)
\psdots[dotsize=1pt 0](27,-11)
\psdots[dotsize=1pt 0](27,-11)
\psdots[dotsize=1pt 0](26,-10)
\psdots[dotsize=1pt 0](26,-10)
\psdots[dotsize=1pt 0](24,-12)
\psdots[dotsize=1pt 0](24,-12)
\psdots[dotsize=1pt 0](25,-13)
\psdots[dotsize=1pt 0](24,-12)
\psdots[dotsize=1pt 0](24,-17)
\psdots[dotsize=1pt 0](25,-16)
\psdots[dotsize=1pt 0](25,-13)
\psdots[dotsize=1pt 0](30,-12)
\psdots[dotsize=1pt 0](29,-13)
\psdots[dotsize=7pt 0](31,-13)
\psdots[dotsize=1pt 0](30,-14)
\psdots[dotsize=1pt 0](30,-12)
\psdots[dotsize=1pt 0](31,-13)
\psdots[dotsize=1pt 0](30,-15)
\psdots[dotsize=1pt 0](31,-16)
\psdots[dotsize=1pt 0](31,-16)
\psdots[dotsize=1pt 0](30,-17)
\psdots[dotsize=1pt 0](30,-17)
\psdots[dotsize=1pt 0](29,-16)
\psdots[dotsize=1pt 0](30,-11)
\psdots[dotsize=1pt 0](30,-12)
\psdots[dotsize=1pt 0](31,-13)
\psdots[dotsize=1pt 0](32,-13)
\psdots[dotsize=1pt 0](31,-16)
\psdots[dotsize=1pt 0](32,-16)
\psdots[dotsize=1pt 0](30,-17)
\psdots[dotsize=1pt 0](30,-18)
\psdots[dotsize=1pt 0](27,-11)
\psdots[dotsize=1pt 0](30,-11)
\psdots[dotsize=1pt 0](30,-11)
\psdots[dotsize=1pt 0](31,-10)
\psdots[dotsize=1pt 0](31,-10)
\psdots[dotsize=1pt 0](33,-12)
\psdots[dotsize=1pt 0](33,-12)
\psdots[dotsize=1pt 0](32,-13)
\psdots[dotsize=1pt 0](26,-10)
\psdots[dotsize=1pt 0](31,-10)
\psdots[dotsize=1pt 0](32,-13)
\psdots[dotsize=1pt 0](32,-16)
\psdots[dotsize=1pt 0](33,-12)
\psdots[dotsize=1pt 0](33,-17)
\psdots[dotsize=1pt 0](33,-17)
\psdots[dotsize=1pt 0](32,-16)
\psdots[dotsize=1pt 0](33,-17)
\psdots[dotsize=1pt 0](31,-19)
\psdots[dotsize=1pt 0](31,-19)
\psdots[dotsize=1pt 0](30,-18)
\psdots[dotsize=1pt 0](30,-18)
\psdots[dotsize=1pt 0](27,-18)
\psdots[dotsize=1pt 0](26,-19)
\psdots[dotsize=1pt 0](31,-19)
\psdots[dotsize=1pt 0](17,-25)
\psdots[dotsize=1pt 0](18,-25)
\psdots[dotsize=1pt 0](18,-25)
\psdots[dotsize=1pt 0](19,-26)
\psdots[dotsize=1pt 0](19,-26)
\psdots[dotsize=1pt 0](19,-27)
\psdots[dotsize=1pt 0](19,-27)
\psdots[dotsize=1pt 0](18,-28)
\psdots[dotsize=1pt 0](18,-28)
\psdots[dotsize=1pt 0](17,-28)
\psdots[dotsize=1pt 0](17,-28)
\psdots[dotsize=1pt 0](16,-27)
\psdots[dotsize=1pt 0](16,-27)
\psdots[dotsize=1pt 0](16,-26)
\psdots[dotsize=1pt 0](16,-26)
\psdots[dotsize=1pt 0](17,-25)
\psdots[dotsize=1pt 0](15,-28)
\psdots[dotsize=1pt 0](16,-29)
\psdots[dotsize=1pt 0](16,-29)
\psdots[dotsize=1pt 0](16,-30)
\psdots[dotsize=1pt 0](15,-28)
\psdots[dotsize=1pt 0](14,-28)
\psdots[dotsize=1pt 0](14,-28)
\psdots[dotsize=1pt 0](13,-29)
\psdots[dotsize=1pt 0](13,-29)
\psdots[dotsize=1pt 0](15,-31)
\psdots[dotsize=1pt 0](15,-31)
\psdots[dotsize=1pt 0](16,-30)
\psdots[dotsize=1pt 0](15,-28)
\psdots[dotsize=1pt 0](16,-27)
\psdots[dotsize=1pt 0](16,-29)
\psdots[dotsize=1pt 0](17,-28)
\psdots[dotsize=1pt 0](16,-26)
\psdots[dotsize=1pt 0](15,-25)
\psdots[dotsize=1pt 0](17,-25)
\psdots[dotsize=1pt 0](16,-23.96)
\psdots[dotsize=1pt 0](16,-23.96)
\psdots[dotsize=1pt 0](15,-25)
\psdots[dotsize=1pt 0](15,-25)
\psdots[dotsize=1pt 0](14,-25)
\psdots[dotsize=1pt 0](16,-23.96)
\psdots[dotsize=7pt 0](16,-23)
\psdots[dotsize=1pt 0](16,-23)
\psdots[dotsize=7pt 0](15,-22)
\psdots[dotsize=1pt 0](15,-22)
\psdots[dotsize=1pt 0](13,-24)
\psdots[dotsize=1pt 0](13,-24)
\psdots[dotsize=1pt 0](14,-25)
\psdots[dotsize=1pt 0](13,-24)
\psdots[dotsize=1pt 0](13,-29)
\psdots[dotsize=1pt 0](14,-28)
\psdots[dotsize=1pt 0](14,-25)
\psdots[dotsize=1pt 0](19,-24)
\psdots[dotsize=1pt 0](18,-25)
\psdots[dotsize=1pt 0](20,-25)
\psdots[dotsize=1pt 0](19,-26)
\psdots[dotsize=1pt 0](19,-24)
\psdots[dotsize=1pt 0](20,-25)
\psdots[dotsize=1pt 0](19,-27)
\psdots[dotsize=1pt 0](20,-28)
\psdots[dotsize=1pt 0](20,-28)
\psdots[dotsize=1pt 0](19,-29)
\psdots[dotsize=1pt 0](19,-29)
\psdots[dotsize=1pt 0](18,-28)
\psdots[dotsize=7pt 0](19,-23)
\psdots[dotsize=1pt 0](19,-24)
\psdots[dotsize=1pt 0](20,-25)
\psdots[dotsize=1pt 0](21,-25)
\psdots[dotsize=1pt 0](20,-28)
\psdots[dotsize=1pt 0](21,-28)
\psdots[dotsize=1pt 0](19,-29)
\psdots[dotsize=1pt 0](19,-30)
\psdots[dotsize=1pt 0](16,-23)
\psdots[dotsize=1pt 0](19,-23)
\psdots[dotsize=1pt 0](19,-23)
\psdots[dotsize=7pt 0](20,-22)
\psdots[dotsize=1pt 0](20,-22)
\psdots[dotsize=1pt 0](22,-24)
\psdots[dotsize=1pt 0](22,-24)
\psdots[dotsize=1pt 0](21,-25)
\psdots[dotsize=1pt 0](15,-22)
\psdots[dotsize=1pt 0](20,-22)
\psdots[dotsize=1pt 0](21,-25)
\psdots[dotsize=1pt 0](21,-28)
\psdots[dotsize=1pt 0](22,-24)
\psdots[dotsize=1pt 0](22,-29)
\psdots[dotsize=1pt 0](22,-29)
\psdots[dotsize=1pt 0](21,-28)
\psdots[dotsize=1pt 0](22,-29)
\psdots[dotsize=1pt 0](20,-31)
\psdots[dotsize=1pt 0](20,-31)
\psdots[dotsize=1pt 0](19,-30)
\psdots[dotsize=1pt 0](19,-30)
\psdots[dotsize=1pt 0](16,-30)
\psdots[dotsize=1pt 0](15,-31)
\psdots[dotsize=1pt 0](20,-31)
\psdots[dotsize=1pt 0](11.51,-7.13)
\psdots[dotsize=1pt 0](16.03,-7.73)
\psdots[dotsize=1pt 0](29.06,-5.24)
\psdots[dotsize=1pt 0](32.97,-5.78)
\psdots[dotsize=1pt 0](6.98,-17.45)
\psdots[dotsize=1pt 0](11.37,-18.06)
\psdots[dotsize=1pt 0](15.22,-14.15)
\psdots[dotsize=1pt 0](26.22,-14.08)
\psdots[dotsize=1pt 0](30.68,-14.69)
\psdots[dotsize=1pt 0](0.17,-7.4)
\psdots[dotsize=1pt 0](4.69,-8)
\psdots[dotsize=1pt 0](20,-14.8)
\psdots[dotsize=1pt 0](12.86,-22.99)
\psdots[dotsize=1pt 0](17.11,-23.6)
\psdots[dotsize=1pt 0](31,-1)
\end{pspicture*}
\end{center}
\caption{Partitioning a Truncated Pyramid (View from Above)}
\label{pyramidtrunc}
\end{figure}

\end{enumerate}

\subsection{Sweeping the $\bcd$-Index}
\label{cdindexsweep}

The partition described in the previous section leads to a recursive
method to compute the $\bcd$-index of $P$ by sweeping.  Each vertex of $P$
will be assigned a certain portion $\Phi_v(P)$ of the $\bcd$-index of
$P$, corresponding to the contribution by $\CB_v(P)$.  This formula is
dual to the results of Stanley~\cite{sta94:fla}.

\begin{thm}
\label{cdsweep}
For any convex $\dd$-polytope $P$,
\begin{enumerate}
\item If $\dd$=0 then $P$ has one vertex $v$ and
$\Phi_v(P)=\Phi(P)=1$.
\item If $\dd>0$ then
\[
\Phi_v(P)=\bd\Phi(R_v)+\sum_{w\in\lver(Q_v)\cap H_v^+}\bc\Phi_w(Q_v),\
v\in\ver(P),
\]
and
\[
\Phi(P)=\sum_{v\in\lver(P)}\Phi_v(P).
\]
\end{enumerate}
\end{thm}

Note in particular that the last vertex $v$ to be swept contributes
nothing to the $\bcd$-index, since $R_v$ is empty, and there are no
vertices $w$ in $\ver(Q_v)\cap H_v^+$.

\

\noindent{\bf Proof.}
We prove by induction that each block in the partition of the faces $T(P)$ has a
$\bcd$-index consisting of a single $\bcd$-word, and that the 
contribution of $\CB_v(P)$ to $\Phi(P)$
is taken into account 
in the formula for $\Phi_v(P)$ stated in the theorem.  This is 
is easy to check for $d=0$: if $P$ is a $0$-polytope with
vertex $v$, then $\CB(P)=\CB_v(P)=\{\{v\}\}$, $\sigma(v)=\emptyset$,
and $\Phi(P)=1$.
So assume $d>0$.

Let $G$ be a middle face as in Step~3 of the partition construction,
and let
$S=\sigma(G)$.  
Note as before that $0\in\sigma(G)$ but $1\not\in\sigma(G)$.
Let $S'=S\setminus\{0\}$.
The four faces that will be in the same pre-block as $G$ will be:
\begin{itemize}
\item $G$, with label set $\{0\}\cup S'$.
\item $\tau(G)$, with label set $\{0,1\}\cup S'$.
\item The face $G'$ for which $\tau(G)$ is the bottom face, with label set $\{1\}\cup S'$.
\item $\tau(G')$, with label set $\{0,1\}\cup S'$.
\end{itemize}
Observe that the label set $\hat S$ of $G\cap H_v$ with respect to the truncation $T(R_v)$ regarded as a $(d-2)$-polytope in its own right is obtained by subtracting $2$ from each label in $S'$. 
Therefore the $\hat S$-chain in $R_v$ contributes in $P$ to 
one $(\{0\}\cup S')$-chain, one $(\{1\}\cup S')$-chain, and two
$(\{0,1\}\cup S')$-chains.  Equation~(\ref{flagh}) then implies that
the contribution to $h_{\{0\}\cup S'}$ and $h_{\{1\}\cup S'}$ is
each $1$.  
Thus, in terms of $\bab$-words, if $u$ is the $\bab$-word for $\hat
S$, then this word contributes $\bb\ba u+\ba\bb u=\bd u$ to the $\bab$-index
of $P$.  
Since such a contribution occurs for each face in a given
block $B$ of $\CB(R_v)$, then the entire block contributes
$\bd\Phi(B)$.  Therefore $\CB(R_v)$ contributes $\bd\Phi(R_v)$ to
$\Phi(P)$.

Now let $G$ be an upper face as in Step~4, and assume $S=\sigma(G)$.  Observe that $0\in\sigma(G)$, and define $S'=S\setminus\{0\}$.  The three faces that will be in the same pre-block as $G$ will be:
\begin{itemize}
\item $G$, with label set $\{0\}\cup S'$.
\item The face $G'$ for which $G$ is the bottom face, with label set $S'$.
\item $\tau(G')$, with label set $\{0\}\cup S'$.
\end{itemize}
Note that the label set $\hat S$ of $G$ with respect to the truncation $T(Q_v)$ regarded as a $(d-1)$-polytope in its own right is obtained by subtracting $1$ from each label in $S'$.
Therefore the $\hat S$-chain in $Q_v$ contributes in $P$ to one $S'$-chain
and two $(\{0\}\cup S')$-chains.  Equation~(\ref{flagh}) then implies
that the contribution to $h_{S'}$ and $h_{\{0\}\cup S'}$ is each $1$.  
Thus, in terms of $\bab$-words, if $u$ is the $\bab$-word for $\hat
S$, then this word contributes $\ba u+\bb u=\bc u$ to the $\bab$-index
of $P$.  
Since such a contribution occurs for each face in a given
block $B$ of $\CB(Q_v)$, then the entire block contributes
$\bc\Phi(B)$.  Therefore $\CB_w(Q_v)$ contributes $\bc\Phi_w(Q_v)$ to
$\Phi(P)$.~$\Box$

\begin{cor}
\label{block}
Each block in the partition of the nonempty faces of $T(P)$ contributes
precisely one $\bcd$-word to $\Phi(P)$.
\end{cor}

\begin{cor}[Stanley]
For a convex $\dd$-polytope $P$ the coefficients of $\Phi(P)$ are
nonnegative.
\end{cor}

\noindent{\bf Examples:}
\begin{enumerate}
\item The line segment ($\dd=1$).  See Figure~\ref{segment}.


If $P$ is a line segment with two
vertices swept in the order $v_1,v_2$, then $Q_{v_i}$ is a point and
$R_{v_i}$ is empty, $i=1,2$.
$Q_{v_1}$ is in $H^+_{v_1}$,
$\Phi_{v_1}(P)=\bc\Phi(Q_{v_1})+\bd\Phi(R_{v_1})=\bc(1)+\bd(0)=\bc$; and
$Q_{v_2}$ is in $H^-_{v_2}$, $\Phi_{v_2}(P)=\bc(0)+\bd(0)=0$.
Thus $\Phi(P)=\bc$.
\item The $n$-gon ($\dd=2$).  See Figure~\ref{polygon}.

\begin{figure}
\begin{center}
\psset{xunit=1.0cm,yunit=1.0cm,algebraic=true,dotstyle=*,dotsize=3pt
0,linewidth=0.8pt,arrowsize=3pt 2,arrowinset=0.25}
\begin{pspicture*}(-1.66,0.08)(5.66,5.66)
\psline(1.36,1.48)(3.32,2.22)
\psline(3.32,2.22)(3.22,4.31)
\psline(3.22,4.31)(1.2,4.87)
\psline(1.2,4.87)(0.05,3.12)
\psline(0.05,3.12)(1.36,1.48)
\psline(1.05,1.86)(1.81,1.65)
\psline(2.88,2.05)(3.3,2.65)
\psline(3.24,3.91)(2.72,4.45)
\psline(1.73,4.72)(0.89,4.39)
\psline(0.34,3.56)(0.32,2.78)
\psline[linestyle=dashed,dash=1pt 1pt](0.22,1.46)(2.3,1.46)
\psline[linestyle=dashed,dash=1pt 1pt](2.18,2.2)(4.26,2.2)
\psline[linestyle=dashed,dash=1pt 1pt](-1.09,3.1)(0.99,3.1)
\psline[linestyle=dashed,dash=1pt 1pt](0.06,4.85)(2.14,4.85)
\psline[linestyle=dashed,dash=1pt 1pt](2.08,4.29)(4.16,4.29)
\rput[tl](1.26,1.22){$v_1$}
\rput[tl](1.74,1.3){$\bc^2$}
\rput[tl](3.24,2.02){$v_2$}
\rput[tl](3.6,2.64){$\bd$}
\rput[tl](-0.3,3){$v_3$}
\rput[tl](-0.26,3.62){$\bd$}
\rput[tl](3.48,4.18){$v_4$}
\rput[tl](3.44,4.84){$\bd$}
\rput[tl](0.86,5.48){$v_5$}
\rput[tl](1.38,5.34){$0$}
\psdots(1.36,1.48)
\psdots(3.32,2.22)
\psdots(3.22,4.31)
\psdots(1.2,4.87)
\psdots(0.05,3.12)
\psdots(1.05,1.86)
\psdots(1.81,1.65)
\psdots(2.88,2.05)
\psdots(3.3,2.65)
\psdots(3.24,3.91)
\psdots(2.72,4.45)
\psdots(1.73,4.72)
\psdots(0.89,4.39)
\psdots(0.34,3.56)
\psdots(0.32,2.78)
\psdots(0.33,3.1)
\psdots(2.87,4.29)
\psdots(2.98,2.2)
\end{pspicture*}
\end{center}
\caption{Sweeping the $\bcd$-Index of a Polygon}
\label{polygon}
\end{figure}
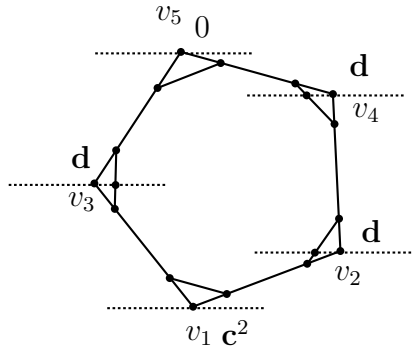

If $P$ is an $n$-gon with vertices swept
in the order $v_1,\ldots,v_n$, then 
$Q_{v_i}$ is a line segment, $i=1,\ldots,n$;
$R_{v_1}$ and $R_{v_n}$ are empty; and $R_{v_i}$ is a point,
$i=2,\ldots,n-1$.
$Q_{v_1}\subset
H^+_{v_1}$, $Q_{v_n}\subset H^-_{v_n}$, and only the top vertex of
$Q_{v_i}$ is in $H^+_{v_i}$, $i=2,\ldots,n-1$.
So $\Phi_{v_1}(P)=\bc\Phi(Q_{v_1})+\bd\Phi(R_{v_1})=\bc(\bc)+\bd(0)=\bc^2$,
$\Phi_{v_n}(P)=\bc(0)+\bd\Phi(R_{v_n})=\bc(0)+\bd(0)=0$, and
$\Phi_{v_i}(P)=\bc(0)+\bd\Phi(R_{v_i})=\bc(0)+\bd(1)=\bd$,
$i=2,\ldots,n-1$.
Thus $\Phi(P)=\bc^2+(n-2)\bd$.
\item 
\label{octahedron}
The octahedron.  

\begin{figure}
\begin{center}
\psset{xunit=1.0cm,yunit=1.0cm,algebraic=true,dotstyle=*,dotsize=3pt
0,linewidth=0.8pt,arrowsize=3pt 2,arrowinset=0.25}
\begin{pspicture*}(-1,-1.68)(6.92,6.68)
\psline[linestyle=dashed,dash=2pt 2pt](1,3)(3,4)
\psline[linestyle=dashed,dash=2pt 2pt](3,4)(4.96,3.24)
\psline(4.96,3.24)(3.16,2.2)
\psline(3.16,2.2)(1,3)
\psline(2,0)(1,3)
\psline[linestyle=dashed,dash=2pt 2pt](3,4)(2,0)
\psline(3.16,2.2)(2,0)
\psline(2,0)(4.96,3.24)
\psline(4,6)(1,3)
\psline[linestyle=dashed,dash=2pt 2pt](4,6)(3,4)
\psline(4,6)(3.16,2.2)
\psline(4,6)(4.96,3.24)
\psline[linestyle=dashed,dash=2pt 2pt](0,-1.52)(5.02,-1.52)
\psline[linestyle=dashed,dash=2pt 2pt](6,-0.54)(5.02,-1.52)
\psline[linestyle=dashed,dash=2pt 2pt](6,-0.54)(0.98,-0.52)
\psline[linestyle=dashed,dash=2pt 2pt](0.98,-0.52)(0,-1.52)
\psline[linestyle=dashed,dash=2pt 2pt]{->}(3,-1)(3,0)
\psline(1.82,0.54)(2.19,0.36)
\psline(2.19,0.36)(2.48,0.52)
\psline[linestyle=dashed,dash=2pt 2pt](2.48,0.52)(2.17,0.68)
\psline[linestyle=dashed,dash=2pt 2pt](2.17,0.68)(1.82,0.54)
\psline(1.16,2.51)(1.44,2.84)
\psline(1.44,2.84)(1.52,3.52)
\psline[linestyle=dashed,dash=2pt 2pt](1.52,3.52)(1.32,3.16)
\psline[linestyle=dashed,dash=2pt 2pt](1.32,3.16)(1.16,2.51)
\psline(2.82,2.33)(2.96,1.83)
\psline(2.96,1.83)(3.48,2.38)
\psline(3.48,2.38)(3.26,2.66)
\psline(3.26,2.66)(2.82,2.33)
\psline(4.48,2.71)(4.37,2.9)
\psline[linestyle=dashed,dash=2pt 2pt](4.78,3.31)(4.76,3.81)
\psline[linestyle=dashed,dash=2pt 2pt](3.28,3.89)(3.14,4.29)
\psline[linestyle=dashed,dash=2pt 2pt](3.14,4.29)(2.71,3.86)
\psline[linestyle=dashed,dash=2pt 2pt](2.71,3.86)(2.88,3.5)
\psline[linestyle=dashed,dash=2pt 2pt](2.88,3.5)(3.28,3.89)
\psline(3.85,5.34)(3.52,5.52)
\psline[linestyle=dashed,dash=2pt 2pt](3.52,5.52)(3.83,5.66)
\psline[linestyle=dashed,dash=2pt 2pt](3.83,5.66)(4.18,5.49)
\psline(4.18,5.49)(3.85,5.34)
\psline(4.37,2.9)(4.76,3.81)
\psline[linestyle=dashed,dash=2pt 2pt](4.78,3.31)(4.48,2.71)
\rput[tl](2.18,-0.06){$v_1$}
\rput[tl](5.3,3.54){$\bcd+\bd\bc$}
\rput[tl](2.3,3.64){$\bd\bc$}
\rput[tl](4.34,6.38){$0$}
\rput[tl](3.32,2.1){$v_2$}
\rput[tl](0.5,3.02){$v_3$}
\rput[tl](5.16,3.08){$v_4$}
\rput[tl](2.64,4.54){$v_5$}
\rput[tl](3.76,6.64){$v_6$}
\rput[tl](0,0.22){$\bc^3+2\bc\bd$}
\rput[tl](3.9,1.92){$2\bcd+\bd\bc$}
\rput[tl](-0.5,3.52){$\bcd+\bd\bc$}
\psdots(1,3)
\psdots(3,4)
\psdots(4.96,3.24)
\psdots(3.16,2.2)
\psdots(2,0)
\psdots(4,6)
\psdots(1.82,0.54)
\psdots(2,0)
\psdots(2.19,0.36)
\psdots(2.48,0.52)
\psdots(2.17,0.68)
\psdots(1.16,2.51)
\psdots(1.44,2.84)
\psdots(1.52,3.52)
\psdots(1,3)
\psdots(1,3)
\psdots(1.32,3.16)
\psdots(2.96,1.83)
\psdots(2.82,2.33)
\psdots(3.48,2.38)
\psdots(3.26,2.66)
\psdots(4.48,2.71)
\psdots(4.37,2.9)
\psdots(4.78,3.31)
\psdots(4.76,3.81)
\psdots(3.28,3.89)
\psdots(3.14,4.29)
\psdots(2.71,3.86)
\psdots(2.88,3.5)
\psdots(3.85,5.34)
\psdots(3.52,5.52)
\psdots(3.83,5.66)
\psdots(4.18,5.49)
\end{pspicture*}
\end{center}
\caption{Sweeping the $\bcd$-Index of an Octahedron}
\label{octahedronfig}
\end{figure}
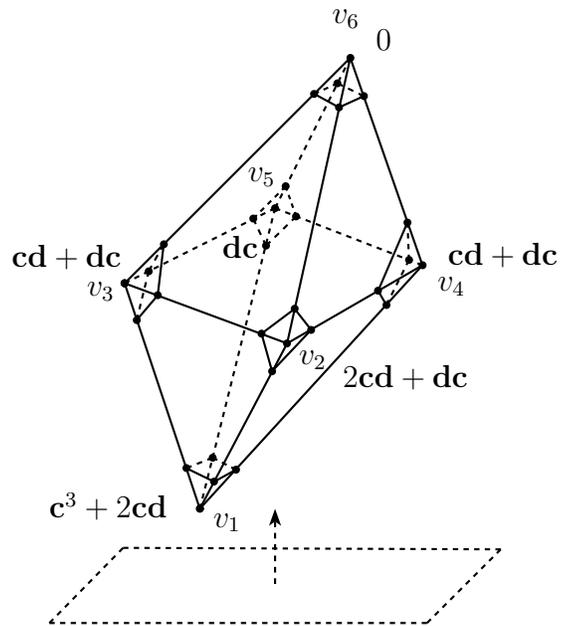

If $P$ is the octahedron with vertices swept in the order
$v_1,\ldots,v_6$ as indicated in Figure~\ref{octahedronfig}, then
$Q_{v_i}$ is a square, $i=1,\ldots,6$; $R_{v_1}$ and $R_{v_6}$ are
empty; and $R_{v_i}$ is a line segment, $i=2,\ldots,5$.
All of the vertices of $Q_{v_1}$ are in $H^+_{v_1}$;
only the top three vertices of $Q_{v_2}$ are in $H^+_{v_2}$;
only the top two vertices of $Q_{v_i}$ are in $H^+_{v_i}$, $i=3,4$;
only the top vertex of $Q_{v_5}$ is in $H^+_{v_5}$;
and none of the vertices of $Q_{v_6}$ are in $H^+_{v_6}$.
So $\Phi_{v_1}(P)=\bc(\bc^2+2\bd)+\bd(0)=\bc^3+2\bcd$, 
$\Phi_{v_2}(P)=\bc(2\bd)+\bd(\bc)=2\bcd+\bd\bc$,
$\Phi_{v_3}(P)=\Phi_{v_4}(P)=\bc(\bd)+\bd(\bc)=\bcd+\bd\bc$,
$\Phi_{v_5}(P)=\bc(0)+\bd(\bc)=\bd\bc$,
and $\Phi_{v_6}(P)=\bc(0)+\bd(0)=0$.
Thus $\Phi(P)=\bc^3+6\bcd+4\bd\bc$ (and we can reverse the letters in
each word of $\Phi(P)$ to get the $\bcd$-index of the cube,
$\bc^3+6\bd\bc+4\bcd$).

\item The square-based pyramid.  See Figure~\ref{pyramidfig2}.

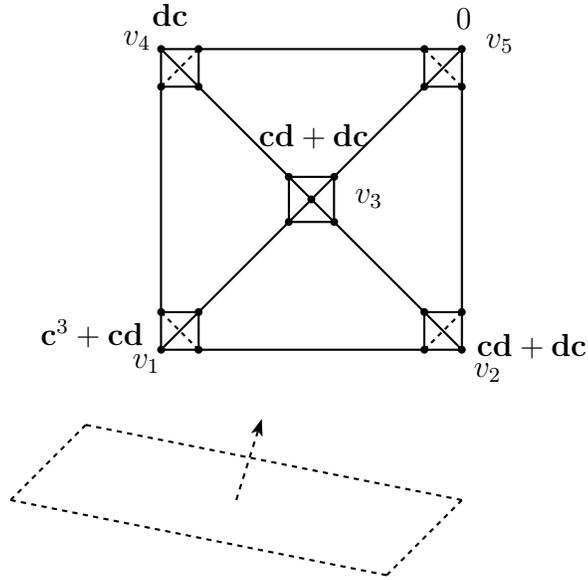
\begin{figure}
\begin{center}
\psset{xunit=1.0cm,yunit=1.0cm,algebraic=true,dotstyle=*,dotsize=3pt
0,linewidth=0.8pt,arrowsize=3pt 2,arrowinset=0.25}
\begin{pspicture*}(-4.38,-5.3)(3.7,2.9)
\psline(-2,2)(2,2)
\psline(2,2)(2,-2)
\psline(2,-2)(-2,-2)
\psline(-2,-2)(-2,2)
\psline(-2,2)(0,0)
\psline(0,0)(2,2)
\psline(2,-2)(0,0)
\psline(0,0)(-2,-2)
\psline(-2,1.5)(-1.5,1.5)
\psline(-1.5,1.5)(-1.5,2)
\psline[linestyle=dashed,dash=2pt 2pt](-2,1.5)(-1.5,2)
\psline(1.5,2)(1.5,1.5)
\psline(1.5,1.5)(2,1.5)
\psline[linestyle=dashed,dash=2pt 2pt](2,1.5)(1.5,2)
\psline(2,-1.5)(1.5,-1.5)
\psline(1.5,-1.5)(1.5,-2)
\psline[linestyle=dashed,dash=2pt 2pt](1.5,-2)(2,-1.5)
\psline(-1.5,-2)(-1.5,-1.5)
\psline(-1.5,-1.5)(-2,-1.5)
\psline[linestyle=dashed,dash=2pt 2pt](-2,-1.5)(-1.5,-2)
\psline(0.3,-0.3)(0.3,0.3)
\psline(0.3,0.3)(-0.3,0.3)
\psline(-0.3,0.3)(-0.3,-0.3)
\psline(-0.3,-0.3)(0.3,-0.3)
\psline[linestyle=dashed,dash=2pt 2pt](-3,-3)(2,-4)
\psline[linestyle=dashed,dash=2pt 2pt](-3,-3)(-4,-4)
\psline[linestyle=dashed,dash=2pt 2pt](-4,-4)(1,-5)
\psline[linestyle=dashed,dash=2pt 2pt](1,-5)(2,-4)
\psline[linestyle=dashed,dash=2pt 2pt]{->}(-1,-4)(-0.66,-2.92)
\rput[tl](-2.38,-2.06){$v_1$}
\rput[tl](1.92,2.54){$0$}
\rput[tl](2.16,-2.18){$v_2$}
\rput[tl](0.58,0.14){$v_3$}
\rput[tl](-2.5,2.24){$v_4$}
\rput[tl](2.32,2.2){$v_5$}
\rput[tl](-3.6,-1.6){$\bc^3+\bcd$}
\rput[tl](2.2,-1.78){$\bcd+\bd\bc$}
\rput[tl](-0.7,1.0){$\bcd+\bd\bc$}
\rput[tl](-2.12,2.6){$\bd\bc$}
\psdots(0,0)
\psdots(-2,2)
\psdots(2,2)
\psdots(2,-2)
\psdots(-2,-2)
\psdots(-2,1.5)
\psdots(-1.5,1.5)
\psdots(-1.5,2)
\psdots(1.5,2)
\psdots(1.5,1.5)
\psdots(2,1.5)
\psdots(2,-1.5)
\psdots(1.5,-1.5)
\psdots(1.5,-2)
\psdots(-1.5,-2)
\psdots(-1.5,-1.5)
\psdots(-2,-1.5)
\psdots(0.3,-0.3)
\psdots(0.3,0.3)
\psdots(-0.3,0.3)
\psdots(-0.3,-0.3)
\end{pspicture*}
\end{center}
\caption{Sweeping the $\bcd$-Index of a Pyramid (View from Above)}
\label{pyramidfig2}
\end{figure}

If $P$ is the square-based pyramid with vertices swept in the order 
$v_1,\ldots,v_5$ as indicated in Figure~\ref{pyramidfig2}, then $Q_{v_i}$ is a
triangle, $i=1,2,4,5$; $Q_{v_3}$ is a square;
$R_{v_1}$ and $R_{v_5}$ are empty; and $R_{v_i}$ is a line segment,
$i=2,3,4$. 
All of the vertices of $Q_{v_1}$ are in $H^+_{v_1}$;
only the top two vertices of $Q_{v_2}$ are in $H^+_{v_2}$;
only the top two vertices of $Q_{v_3}$ are in $H^+_{v_3}$;
only the top vertex of $Q_{v_4}$ is in $H^+_{v_4}$;
and none of the vertices of $Q_{v_5}$ are in $H^+_{v_5}$.
So $\Phi_{v_1}(P)=\bc(\bc^2+\bd)+\bd(0)=\bc^3+\bcd$, 
$\Phi_{v_2}(P)=\bc(\bd)+\bd(\bc)=\bcd+\bd\bc$,
$\Phi_{v_3}(P)=\bc(\bd)+\bd(\bc)=\bcd+\bd\bc$,
$\Phi_{v_4}(P)=\bc(0)+\bd(\bc)=\bd\bc$,
and $\Phi_{v_5}(P)=\bc(0)+\bd(0)=0$.
Thus $\Phi(P)=\bc^3+3\bcd+3\bd\bc$.

\end{enumerate}

\subsection{A Symmetric Formula}

Since the $\bcd$-index is independent of the sweeping used, we can
symmetrize the formula in Theorem~\ref{cdsweep}
by taking the average of the results from a
sweep and its opposite.  
In the following theorem the contribution $\Phi_v(P)$ from the sweep 
is different from that in Theorem~\ref{cdsweep}, even though we are using 
the same notation.  Note in particular that $\Phi_v(P)$ now involves
the entire $\bcd$-indices of both $Q_v$ and $R_v$.

\begin{thm}
\label{cdsymm}
For any convex $\dd$-polytope $P$,
\begin{enumerate}
\item If $\dd$=0 then $P$ has one vertex $v$ and
$\Phi_v(P)=\Phi(P)=1$.
\item If $\dd>0$ then
\[
\Phi_v(P)=\frac{1}{2}[ \bc\Phi(Q_v)+(2\bd-\bc^2)\Phi(R_v)],\ v\in\ver(P),
\]
and
\[
\Phi(P)=\sum_{v\in\lver(P)}\Phi_v(P).
\]
\end{enumerate}
\end{thm}

\noindent{\bf Proof.}
It is helpful first to extend the computation of the $\bcd$-index to 
some ``near'' polytopes.  
Let $R$ be a $(d-1)$-polytope and consider
the infinite cylinder $R\times\BR$ with two points $v^+$ and $v^-$ adjoined at
infinity, one in each direction, each declared to be
formally incident to each of the faces of
the cylinder.  Call this object $\oR$.  Now $\oR$ is not a $d$-polytope,
but its complete truncation $T(\oR)$ is:  first truncate each of its two
vertices by capping the cylinder with a hyperplane at each end,
resulting in a prism over $R$.  Then continue by truncating the faces
of $R$.  In sweeping the $\bcd$-index of $R$ from $v^-$ toward $v^+$,
the last vertex $v^+$ contributes nothing.  Now $R_{v^-}$ is empty and
$Q_{v^-}$ is combinatorially equivalent to the original $R$.
Therefore by Theorem~\ref{cdsweep}, $\Phi(\oR)=\bc\Phi(R)$.

Now let $P$ be a $d$-polytope with vertices swept in the order
$v_1,\ldots,v_\ell$.  
For each vertex $v$ 
define $\rPhi_v(P)$ to be the contribution by $v$ to $\Phi(P)$ in this
sweeping order, and $\lPhi_v(P)$ to be the contribution by $v$ to the
$\bcd$-index of $P$ in the reverse sweeping direction.
Hence
\[
\Phi(P)=\sum_{i=1}^\ell\rPhi_{v_i}(P)=\sum_{i=1}^\ell\lPhi_{v_i}(P).
\]

Let $H$ be a hyperplane in the sweeping family
positioned so that it separates $v_k$ from $v_{k+1}$.
Define $P^+$ to be the object obtained by taking $P\cap H^-$, 
applying
a projective transformation that sends the facet $P\cap H$ to
infinity, and adjoining a point $v^+$ at infinity, formally incident to each of
the unbounded faces of $P^+$.  (This latter operation is dual to the
``capping'' operation arising in $S$-shellings\@.)
Again $P^+$ is not a polytope, but its complete truncation $T(P^+)$ is:
first truncate $v^+$ by capping the unbounded faces of $P^+$ with a
single hyperplane.  Then continue by truncating the other vertices, and then
the other faces.
In sweeping the $\bcd$-index of $P^+$ in the same vertex order as $P$,
the last vertex $v^+$ contributes nothing, and the remaining vertices
contribute to the $\bcd$-index of $P^+$ in the same way that they
contributed to $P$.  Thus
\[
\Phi(P^+)=\sum_{i=1}^k\rPhi_{v_i}(P).
\]
In a similar manner, define $P^-$ by taking $P\cap H^+$, 
applying
a projective transformation that sends the facet $P\cap H$ to
infinity, and adjoining a point $v^-$ at infinity, formally incident to each of
the unbounded faces of $P^-$.  
Then
\[
\Phi(P^-)=\sum_{i=k+1}^\ell\lPhi_{v_i}(P).
\]

Let $R=P\cap H$.
Now as complexes, $P^+$ and $P^-$ together equal $P$ with an
extra copy of $\oR$, so
\[
\Phi(P^+)+\Phi(P^-)=\Phi(P)+\Phi(\oR)=\Phi(P)+\bc\Phi(R).
\]
Thus
\[
\begin{array}{rcl}
\Sum_{i=k+1}^\ell\rPhi_{v_i}(P)+\Sum_{i=1}^k\lPhi_{v_i}(P)&=&
2\Phi(P)-\Sum_{i=1}^k\rPhi_{v_i}(P)-\Sum_{i=k+1}^\ell\lPhi_{v_i}(P)\\\ \\
&=&2\Phi(P)-(\Phi(P)+\bc\Phi(R))\\\ \\
&=&\Phi(P)-\bc\Phi(R).
\end{array}
\]
Applying the above formula to $Q_v$,
Theorem~\ref{cdsweep}
then implies
\[
\begin{array}{rcl}
\Phi_v(P)&=&\frac{1}{2}[\Phi_v(P)+\Phi_v(P)]\\\ \\
&=&\frac{1}{2}
[\bd\Phi(R_v)+\Sum_{w\in\lver(Q_v)\cap H_v^+}\bc\rPhi_w(Q_v)
+\bd\Phi(R_v)+\Sum_{w\in\lver(Q_v)\cap H_v^-}\bc\lPhi_w(Q_v)]\\\ \\
&=&\frac{1}{2}[2\bd\Phi(R_v)+\bc\Phi(Q_v)-\bc^2\Phi(R_v)]\\\ \\
&=&\frac{1}{2}[\bc\Phi(Q_v) + (2\bd-\bc^2)\Phi(R_v) ].~\Box
\end{array}
\]
Though it might not be obvious from the formula,
note that $\Phi_v(P)$ in the theorem
is necessarily nonnegative since it is the sum of two nonnegative
quantities.

\ 

\noindent{\bf Examples:}
\begin{enumerate}
\item The line segment.  See Figure~\ref{segment}.
$\Phi_{v_i}(P)=\frac{1}{2}[\bc\Phi(Q_{v_i})+(2\bd-\bc^2)\Phi(R_{v_i})]
=\frac{1}{2}[\bc(1)+(2\bd-\bc^2)(0)]=\frac{1}{2}\bc$, $i=1,2$.
Thus $\Phi(P)=\bc$.
\item The $n$-gon.  See Figure~\ref{polygon}.
For $i=1$ or $i=n$,
$\Phi_{v_i}(P)=\frac{1}{2}[\bc\Phi(Q_{v_i})+(2\bd-\bc^2)\Phi(R_{v_i})]=
\frac{1}{2}[\bc(\bc)+(2\bd-\bc^2)(0)]=\frac{1}{2}\bc^2$;
and for $i=2,\ldots,n-1$,
$\Phi_{v_i}(P)=\frac{1}{2}[\bc\Phi(Q_{v_i})+(2\bd-\bc^2)\Phi(R_{v_i})]=\frac{1}{2}[\bc^2+(2\bd-\bc^2)]=\bd$,
$i=2,\ldots,n-1$.
Thus $\Phi(P)=\bc^2+(n-2)\bd$.
\item 
The octahedron.  See Figure~\ref{octahedron}.
$\Phi_{v_i}(P)=\frac{1}{2}[\bc(\bc^2+2\bd)+(2\bd-\bc^2)(0)]=\frac{1}{2}\bc^3+\bcd$,
$i=1$ and $i=6$;
and
$\Phi_{v_i}(P)=\frac{1}{2}[\bc(\bc^2+2\bd)+(2\bd-\bc^2)(\bc)]=\bcd+\bd\bc$, $i=2,\ldots,5$.
Thus $\Phi(P)=\bc^3+6\bcd+4\bd\bc$.

\item The square-based pyramid.  See Figure~\ref{pyramidfig}.
$\Phi_{v_i}(P)=\frac{1}{2}[\bc(\bc^2+\bd)+(2\bd-\bc^2)(0)]=\frac{1}{2}\bc^3+\frac{1}{2}\bcd$,
$i=1$ and $i=5$;
$\Phi_{v_i}(P)=\frac{1}{2}[\bc(\bc^2+\bd)+(2\bd-\bc^2)(\bc)]=\frac{1}{2}\bcd+\bd\bc$,
$i=2,4$;
and
$\Phi_{v_3}(P)=\frac{1}{2}[\bc(\bc^2+2\bd)+(2\bd-\bc^2)(\bc)]=\bcd+\bd\bc$,
Thus $\Phi(P)=\bc^3+3\bcd+3\bd\bc$.
\end{enumerate}

\section{The Toric $h$-Vector}

\subsection{Definitions}

The {\em toric $h$-vector\/} 
of (the boundary complex of) a convex $\dd$-polytope $P$, 
$h(\partial P)=(h_0,\ldots,h_\dd)$, is a linear combination of the
components of the flag $h$-vector that is a
nonnegative, symmetric,
generalization of the $h$-vector of a simplicial polytope.  The
component $h_i=h_i(\partial P)$ is the rank of the $(2\dd-2i)$th
middle perversity intersection homology group of the associated toric
variety in the case that $P$ is
rational (has a realization with rational vertices).  The
$g$-Theorem~\cite{sta80:num} implies that the $h$-vector
of a simplicial polytope is unimodal.
Karu~\cite{kar04:har} 
proved that this is also the case for the toric $h$-vector of a
general polytope $P$, even when $P$ is not rational.  For a summary of
some other results on the toric $h$-vector see~\cite{baylee93:com}.

To define the toric $h$-vector recursively, let
$h(\partial P,x)=\sum_{i=0}^\dd h_ix^{\dd-i}$ and 
$g(\partial P,x)=\sum_{i=0}^{\lfloor
\dd/2\rfloor}g_ix^i$ where $g_0=g_0(\partial P)=h_0$ and
$g_i=g_i(\partial P)=h_i-h_{i-1}$,
$i=1,\ldots,\lfloor \dd/2\rfloor$.  Then
\[
g(\emptyset,x)=h(\emptyset,x)=1,
\]
and
\[
h(\partial P,x)=\sum_{{\mbox{\scriptsize $G$ face of $\partial P$}}}
       g(\partial G,x)(x-1)^{\dd-1-\dim G}.
\]
In the case that $P$ is simplicial the toric $h$-vector of
$\partial P$ agrees with the simplicial $h$-vector of $P$.

For example, the toric $h$-vectors of the boundary complexes of a 
point, line segment, $n$-gon, 
octahedron, and cube are, respectively, $(1)$, $(1,1)$, 
$(1,n-2,1)$, $(1,3,3,1)$, and $(1,5,5,1)$. 

\subsection{Sweeping the Toric $h$-Vector}

In Section~\ref{ordinaryh} we recalled that by sweeping any simple polytope
$P$ we can compute the $h$-vector of its dual $P^*$.
Analogously, as we sweep any polytope $P$, we can compute the toric 
$h$-vector of its dual $P^*$.

Define operators $\bc:\BR^{d+1}\rightarrow\BR^{d+2}$ and 
$\bd:\BR^{\dd+1}\rightarrow\BR^{\dd+3}$ on symmetric vectors
$(h_0,\ldots,h_d)$ by
\[
(h_0,\ldots,h_\dd)\bc=\left\{
\begin{array}{ll}
(g_0,g_1,\ldots,
g_{\lfloor \dd/2\rfloor},
g_{\lfloor \dd/2\rfloor},
\ldots,
g_1,g_0)
&
\mbox{if $\dd$ is even}\\
(g_0,g_1,\ldots,
g_{\lfloor \dd/2\rfloor},
0,
g_{\lfloor \dd/2\rfloor},
\ldots,
g_1,g_0)
&
\mbox{if $\dd$ is odd}
\end{array}
\right.
\]
and
\[
(h_0,\ldots,h_\dd)\bd=\left\{
\begin{array}{ll}
(0,\ldots,0,g_{\lfloor \dd/2\rfloor},0,\ldots,0)
&\mbox{if $\dd$ is even}\\
(0,\ldots,0)
&\mbox{if $\dd$ is odd}
\end{array}
\right.
\]
where as before $g_0=h_0$ and $g_i=h_i-h_{i-1}$, $i=1,\ldots,\lfloor
\dd/2\rfloor$.

Define (with a small abuse of notation)
$h_v(\partial P^*)$ to be the contribution by $v$ to the toric
$h$-vector of $P^*$ during the sweeping of $P$.
We now
have an analog to Theorem~\ref{cdsweep}:

\begin{thm}
\label{cdtoh0}
For any convex $\dd$-polytope $P$,
\begin{enumerate}
\item If $\dd=0$ then $P$ has one vertex $v$ and
$h_v(\partial P^*)=h(\partial P^*)=(1)$.
\item If $\dd>0$ then, regarding $\bc$ and $\bd$ as operators,
\[
h_v(\partial P^*)=h(\partial (R_v)^*)\bd+\sum_{w\in\lver(Q_v)\cap
H^+_v}h_w(\partial (Q_v)^*)\bc,\
v\in\ver(P),
\]
and
\[
h(\partial P^*)=\sum_{v\in\lver(P)}h_v(\partial P^*).
\]
\end{enumerate}
\end{thm}

\noindent{\bf Proof.}
Returning to the definitions of the operators $\bc$ and $\bd$, write
$h(x)=\sum_{i=0}^dh_ix^i$ and
$g(x)=\sum_{i=0}^{\lfloor\frac{d}{2}\rfloor}g_ix^i$.  
For any polynomial $p(x)=\sum_{i=0}^dp_ix^i$ and nonnegative integer $m$
define $U_{\leq m}p(x)=\sum_{i=0}^mp_ix^i$.
Then it is easy to verify that the operators $\bc$ and $\bd$ can be expressed as
\[
\begin{array}{rcl}
h(x)\bc&=&(x-1)h(x)+2g(x),\\
h(x)\bd&=&(x-1)g(x)+U_{\leq m}[(1-x)g(x)],
\end{array}
\]
where $m=\lfloor\frac{d+1}{2}\rfloor$.
Bayer and Ehrenborg~\cite{bayehr00:tor} developed explicit formulas for
computing the toric $h$-vector from the $\bcd$-index (Theorem~4.2) in
which the contribution for each $\bcd$-word is determined.  Their
Lemma~7.9 and Proposition~7.10 relate the contribution toward the
toric $h$-vector for
$\bcd$-words $u\bc$ and $u\bd$ with that of $\bcd$-word $u$, and these
correspond precisely to the formulas for the
operators $\bc$ and $\bd$ defined above.

For any $\bcd$-polynomial $\Phi$ write $\Phi^*$ for the polynomial
resulting from reversing all of the words in $\Phi$.  Thus for any
polytope $P$, $\Phi(P^*)=\Phi^*(P)$.

By Theorem~\ref{cdsweep},
\[
\Phi(P^*)=\Phi^*(P)=\sum_{v\in\lver(P)}\Phi^*_v(P),
\]
and
\[
\Phi^*_v(P)=\Phi^*(R_v)\bd+\sum_{w\in\lver(Q_v)\cap
H_v^+}\Phi^*_w(Q_v)\bc,\
v\in\ver(P).
\]
Now use induction and compute the toric $h$-vectors of both sides.~$\Box$

Induction immediately yields a formula to obtain the toric
$h$-vector directly from the $\bcd$-index
and to an analog of
Theorem~\ref{cdsymm}.

\begin{thm}
\label{cdtoh1}
Let $P$ be a convex $\dd$-polytope.  Then, regarding $\bc$ and $\bd$ 
as operators,
$h(\partial P)=(1)\Phi(P)$.
\end{thm}

Lemma~7.9 and Proposition~7.10 of \cite{bayehr00:tor} can be regarded
as definitions of operators 
$\bc$ and $\bd$ 
acting upon toric $h$-vectors, and these results imply
Theorem~\ref{cdtoh1} directly.

In the following theorem the contribution $h_v(\partial P^*)$ from the sweep
is different from that in Theorem~\ref{cdtoh0}, even though we are
using
the same notation.  Note in particular that $h_v(\partial P^*)$ now involves
the entire toric $h$-vectors of both $\partial(Q_v)^*$ and
$\partial(R_v)^*$.

\begin{thm}
\label{cdtoh2}
\label{hsymm}
For any convex $\dd$-polytope $P$,
\begin{enumerate}
\item If $\dd=0$ then $P$ has one vertex $v$ and $h_v(\partial P^*)=h(\partial P^*)=(1)$.
\item If $\dd>0$ then, regarding $\bc$ and $\bd$ as operators,
\[
h_v(\partial P^*)=\frac{1}{2}[h(\partial (Q_v)^*)\bc+h(\partial
(R_v)^*)(2\bd-\bc^2)],\ v\in\ver(P),
\]
and
\[
h(\partial P^*)=\sum_{v\in\lver(P)}h_v(\partial P^*).
\]
\end{enumerate}
\end{thm}

\noindent{\bf Examples}
\begin{enumerate}
\item If $\dd=0$ and $P$ is a point then $h(\partial P)=(1)\Phi(P)=(1)1=(1)$.
\item If $\dd=1$ and $P$ is a line segment then $h(\partial P)=(1)\bc=(1,1)$.
\item If $\dd=2$ and $P$ is an $n$-gon then
\[
\begin{array}{rcl}
h(P)&=&(1)\Phi(P)\\
&=&(1)(\bc^2+(n-2)\bd)\\
&=&(1,1)\bc+(n-2)(0,1,0)\\
&=&(1,0,1)+(n-2)(0,1,0)\\
&=&(1,n-2,1).
\end{array}
\]
We can also use Theorem~\ref{cdtoh0}; see Figure~\ref{polygon}.
Vertex $v_1$ contributes $(1,1)\bc=(1,0,1)$ and each remaining vertex
except the last contributes $(1)\bd=(0,1,0)$, yielding $(1,n-2,1)$.
\item If $\dd=3$ and $P$ is the cube then
\[
\begin{array}{rcl}
h(\partial P^*)&=&\Phi(P)(1)\\
&=&(\bc^3+6\bcd+4\bd\bc)(1)\\
&=&\bc^2(1,1)+6\bc(0,1,0)+4\bd(1,1)\\
&=&\bc(1,0,1)+6(0,1,1,0)+4(0,0,0,0)\\
&=&(1,-1,-1,1)+(0,6,6,0)+(0,0,0,0)\\
&=&(1,5,5,1).
\end{array}
\]

We can use Theorem~\ref{cdtoh0} to compute the toric $h$-vector of a cube
$P^*$ from a sweeping of the octahedron $P$
(see Figure~\ref{octahedron}):
$h_{v_1}(\partial P^*)=(1,2,1)\bc+(0)\bd=(1,1,1,1)$, 
$h_{v_2}(\partial P^*)=(0,2,0)\bc+(1,1)\bd=(0,2,2,0)+(0,0,0,0)=(0,2,2,0)$,
$\Phi_{v_3}(P)=\Phi_{v_4}(P)=(0,1,0)\bc+(1,1)\bd=(0,1,1,0)+(0,0,0,0)=(0,1,1,0)$,
$\Phi_{v_5}(P)=(0)\bc+(0,1,1,0)\bd=(0,0,0,0)$,
and $\Phi_{v_6}(P)=(0)\bc+(0)\bd=0$.
Thus $h(\partial P^*)=(1,5,5,1)$.

We can also apply Theorem~\ref{hsymm} to the octahedron to compute the
$h$-vector of the cube:
$h_{v_i}(\partial P^*)=
\frac{1}{2}[(1,2,1)\bc+(0,0)(2\bd-\bc^2)]=\frac{1}{2}(1,1,1,1)$, $i=1$
and $i=6$; and
$h_{v_i}(P^*)=\frac{1}{2}[(1,2,1)+\bc(1,1)(2\bd-\bc^2)]=\frac{1}{2}[(1,1,1,1)
+2(0,0,0,0)-(1,-1,-1,1))]=\frac{1}{2}(0,2,2,0)=(0,1,1,0)$,
$i=2,\ldots,5$.  Thus $h(\partial P^*)=(1,5,5,1)$.
\end{enumerate}

\subsection{An ``Extended Toric'' $h$-Vector}
\label{exttoric}

Even though for a $\dd$-polytope $P$ the $\bcd$-index $\Phi(P)$ contains
$F_\dd-1$
independent pieces of information, the toric $h$-vector $h(P)$
contains only $\lfloor(\dd+1)/2\rfloor$ independent pieces of information.
The source of the loss from $\Phi(P)$ to $h(P)$ is evident---the $\bd$
operator ``erases'' information.  We can get around this by keeping
track of some of the intermediate calculations (those vectors that are
about to be acted upon by $\bd$).

Let $W$ be the set of all $\bcd$-words $w$ of degree at most $\dd$
(including the word $1$).  Denote by $W^{\lbd}$ the set of all words in $W$
having $\bd$ as the first letter, and include $1$ in this set also.
For $w\in W$ let $\Phi^w(P)w$ be that portion of $\Phi(P)$ with terms
ending in $w$.  Define $h^w(P)=(1)\Phi^w(P)$.  
Define the {\em ``extended toric'' $h$-vector\/} of $P$ to be
$\hh(P)=(h^w(P):w\in W^\lbd)$.

For example, if $P$ is the octahedron, then 
$\Phi(P)=\bc^3+4\bd\bc+6\bcd$.  We have:
\[
\begin{array}{ccc}
w&\Phi^w(P)&h^w(P)\\\hline
1&\bc^3+4\bd\bc+6\bcd&(1,3,3,1)\\
\bc&\bc^2+4\bd&(1,4,1)\\
\bd&6\bc&(6,6)\\
\bc^2&\bc&(1,1)\\
\bd\bc&4&(4)\\
\bcd&6&(6)\\
\bc^3&1&(1)
\end{array}
\]
Then $W^\lbd=\{1, \bd, \bd\bc\}$ and the extended toric $h$-vector is
$\hh(P)=(h^1(P),h^\lbd(P),h^{\lbd\lbc}(P))=((1,3,3,1),(6,6),(4)))$.

\begin{thm}
For a $\dd$-polytope $P$ each $h^w(P)$, $w\in W^\lbd$,
is nonnegative, symmetric, and
unimodal, and $\hh(P)$ determines $\Phi(P)$.
\end{thm}

To prove this, recall that the toric $h$-vector of any polytope is
nonnegative, symmetric, and unimodal, and by the recursive application
of Proposition~\ref{cdsweep}
the operator $\bd$ is always multiplied onto the $\bcd$-index of some
polytope.  Hence each $h^w(P)$, $w\in W^\lbd$, 
being a sum of $h$-vectors of such
polytopes, is nonnegative, symmetric, and unimodal.  
To show that $\hh(P)$
determines $\Phi(P)$, observe that
\begin{enumerate}
\item Any symmetric vector $h$ can be recovered from $h\bc$.
\item For any $\bcd$-word $w$, $h^{\lbc w}(P)$ can be recovered from
$h^w(P)$ and $h^{\lbd w}(P)$, since 
$h^w(P)=(h^{\lbc w}(P))\bc+(h^{\lbd w}(P))\bd$.
Therefore, by reverse induction on the degree of $w$, we can
recover all of the vectors $h^w(P)$ from $\hh(P)$.
\item For any $\bcd$-word $w$ of degree $\dd$, the coefficient of $w$ in
$\Phi(P)$ is precisely the single entry of $h^w(P)$.
\end{enumerate}
This concludes the proof.~$\Box$

At this point it remains to be seen whether or not one can get a
better understanding of the collection of flag $f$-vectors of convex
$\dd$-polytopes from their extended toric $h$-vectors, or indeed whether
one is even justified in giving $\hh(P)$ this name.  

\section{Comments}

Karu~\cite{kar06:cd} described the $\bcd$-index of a complete fan
$\Delta$ by
beginning with its first barycentric subdivision which, in the case of
polytopes, is dual to the complete truncation.  He defines operators
$C$ and $D$ on functions $f:\Delta^{\leq m}\rightarrow{\BZ}$ on the
$m$-skeleta of the fan $\Delta$.  He proves (Proposition~1.2) that
if $u$ is a $\bcd$-word, then the result of applying the corresponding
$CD$
operator to the constant function $1$ on $\Delta$ is the
coefficient of $u$ in the $\bcd$-index of $\Delta$.  He then
demonstrates how $C$ and $D$ have counterparts in the category of
sheaves, and uses this to prove nonnegativity of the $\bcd$-index of
$\Delta$.
Karu asks what the coefficients of the $\bcd$-index count, and so we
can now provide one answer of a sort
in the case of complete fans associated
with polytopes.
It is natural to ask what the connection might be between the operators $C$
and $D$ and the toric $h$-vector.


\end{sloppypar}
\end{document}